\magnification = 1200
\input amssym.def
\input amssym.tex
\input graphicx
\def \qed {\hfill $\square$}
\def \R {\Bbb R}

\def \d {\partial}
\def \H {{\cal H}}
\def \ms {\medskip}
\def \msi {\medskip\noindent}
\def \ssi {\smallskip\noindent}
\def \ss {\smallskip}

\def \pari {\par\noindent}
\def \vfe {\vfill\eject}

\def \sm {\setminus }

\def \dsp {\displaystyle }
\def \wt {\widetilde }

\overfullrule=0pt
\def \dist{\mathop{\rm dist}\nolimits}
\def \diam{\mathop{\rm diam}\nolimits}

\def \inf{\mathop{\rm inf}\nolimits}

\centerline{\bf SHOULD WE SOLVE PLATEAU'S PROBLEM AGAIN?}
\smallskip
\centerline{Guy David}

\vskip 1cm
\noindent
{\bf R\'{e}sum\'{e}.}
Apr\`{e}s une courte description de plusieurs solutions 
classiques du probl\`{e}me de Plateau, on parle d'autres 
mod\'{e}lisations des films de savon, et de probl\`{e}mes
ouverts li\'{e}s. On insiste un peu plus sur un mod\`{e}le
bas\'{e} sur des d\'{e}formations et des conditions glissantes
\`{a} la fronti\`{e}re.

\bigskip \noindent
{\bf Abstract.}
After a short description of various classical solutions 
of Plateau's problem, we discuss other ways to model soap films, 
and some of the related questions that are left open. 
A little more attention is payed to a more specific model, with 
deformations and sliding boundary conditions.


\ms\ms
\msi{\bf 1. Introduction}
\ss
The main goal of this text is to give a partial account
of the situation of Plateau's problem, on the existence and
regularity of soap films with a given boundary. 
We intend to convince the reader that there are many 
reasonable ways to state a Plateau problem, most of which
give interesting questions that are still wide open.
This is even more true when we want our models to stay
close to Plateau's original motivation, which was to describe
physical phenomena such as soap films.

Plateau problems led to lots of beautiful results;
we shall start the paper with a rapid description
of some of the most celebrated solutions of Plateau's 
problem (Section~2), followed by a description of a 
few easy examples (Section 3), mostly to explain more
visibly some objections and differences between the models.

With these examples in mind, we shall shortly return to 
the modeling problem, and mention a few additional
ways to state a Plateau problem and (in some cases)
get solutions with a nice physical flavor (Section 4).

It turns out that the author has a preference for a specific
way to state Plateau problems, coming from Almgren's
notion of minimal sets, which we accommodate at the boundary
with what we shall call sliding boundary conditions. 

In Section 5, we shall describe briefly the known
local regularity properties of the Almgren minimal sets
(i.e., far from the boundary), and why we would like to extend
some of these regularity results to sliding minimal sets,
all the way to the boundary. We want to do this because little
seems to be known on the boundary behavior of soap films
and similar objects, but also because we hope that this 
may help us get existence results.

We try to explain this in Section 6, and at the same time why, 
even though beautiful compactness and stability results 
for various classes of objects (think about Almgren minimal sets, 
but also currents or varifold) yield relatively systematic
partial solutions of Plateau problems, these solutions are 
not always entirely satisfactory (we shall call these amnesic solutions).

We explain in Section 7 why the regularity results for
sliding Almgren minimal sets also apply to solutions of the
Reifenberg and size minimization problems described in Section 2.

The author wishes to thank K. Brakke and John M. Sullivan for allowing 
him to use beautiful pictures from their site, T. De Pauw, J. Harrison, 
and F. Morgan for help in the preparation of this manuscript, and 
the organizers of the Stein conference in Princeton for a wonderful
event. Special congratulations and thanks are due to E. Stein who is 
an inspiring example as a mathematician and the leader of an amazing school.

\msi{\bf 2. SOME CELEBRATED SOLUTIONS}
\ss
In this section we describe some of the most celebrated 
ways to state a Plateau problem and solve it.
To make things easier, we shall mostly think about the simplest
situation where we give ourselves a smooth simple loop $\Gamma$ in
$\R^3$, and we look for a surface bounded by $\Gamma$, with
minimal area. Even that way, a few different definitions of the terms
``bounded by" and ``area" will be used.

\ssi{\bf 2.a. Parameterizations by disks, Garnier, Douglas, and 
Rad\'{o}} 
\ss
Here we think of surfaces as being parameterized, and 
compute their area as the integral of a Jacobian determinant.
For instance, let $\Gamma \i \R^n$ be a simple curve, 
which we parameterize with a continuous function 
$g: \d D \to \Gamma$, where $D$ denotes the unit disk 
in $\R^2$; we decide to minimize the area
$$
A(f) = \int_{D} J_f(x) dx,
\leqno (2.1)
$$
where $J_f(x)$ denotes the positive Jacobian of $f$ at $x$,
among a suitable class of functions 
$f : \overline D \to \R^n$ such that $f_{|\d B} = g$.
In fact, it is probably a good idea to allow also
functions $f$ such that $f_{|\d B}$ is another equivalent 
parameterization of $\Gamma$. 

There are obvious complications with this problem,
and the main one is probably the lack of
compactness of the reasonable classes of acceptable 
parameterizations. That is, if $\{ f_k \}$ is 
a minimizing sequence, i.e., if $A(f_k)$ tends
to the infimum of the problem, and even if the 
sets $f_k({\Bbb D})$ converge very nicely to a beautiful
smooth surface, the parameterizations $f_k$ themselves
could have no limit. Even if $\Gamma$ is the unit circle
and each $f_k(\Bbb D)$ is equal to $D$, 
it could be that we stupidly took a sequence of smooth 
diffeomorphisms $f_k : D \to D$ that behave more and 
more wildly.

For $2$-dimensional surfaces, there is at least
one standard way to deal with this problem: we can
decide to use conformal parameterizations of the image,
normalized in some way, and gain compactness this way.
This is more or less the approach that was taken,
for instance, by Garnier [Ga], 
and then Tibor Rad\'{o} [Ra1] (1930) 
and J. Douglas [Do] (1931). 

Let us say a few words about the existence
theorem of Douglas, who was able to prove the
existence of a function $f$ that minimizes $A(f)$
under optimal regularity conditions for $\Gamma$.

Let us say a few words about the idea because it 
is beautiful. We decide that $f$ will be the harmonic 
extension of its boundary values on $\d D$, which we 
require to be a parameterization of $\Gamma$ 
(but not given in advance). This is a reasonable thing to do, 
because we know that such harmonic parameterizations exist, 
at least in the smooth case. Then we can compute $A(f)$
in terms of $g = f_{|\d B}$, and we get that
$A(f) = B(g)$, where
$$
B(g) = \int_0^{2\pi}\int_0^{2\pi} 
\ {\sum_{j=1}^n |g_j(\theta)-g_j(\varphi)|^2
\over \sin^2\big({\theta-\varphi \over 2}\big)} 
\, d\theta d\varphi,
\leqno (2.2)
$$
where the $g_j$ are just the coordinates of $g$.

Now $B$ is a much easier functional to minimize,
in particular because there is one less variable,
and Douglas obtains a solution rather easily. 
The paper [Do] seems very simple and pleasant to read.

Of course the mapping $f$ is smooth, but there were 
still important regularity issues to be resolved, 
concerning the way $f(D)$ is embedded, or whether
$f$ may have critical points. See
[Laws], [Ni], or [Os].  
 
There are two or three important difficulties with this way
of stating Plateau's problem. The minor one is that
getting reasonably normalized parameterizations 
will be much harder for higher dimensional sets,
thus making existence results in these dimensions
much less likely.

Even when the boundary is a nice curve $\Gamma$, many of the
physical solutions of Plateau's problem are not parameterized by a 
disk, but by a more complicated set, typically a Riemann surface 
with a boundary. So we should also allow more domains than 
just disks; this is not such a serious issue though.

But also, in many cases the solutions of Douglas do not
really describe soap films (which were at the center
of Plateau's initial motivation). For instance, 
if $\Gamma$ is folded in the right way, the minimizing 
surface $f(D)$ of Douglas will cross itself, 
and because we just minimize the integral of the Jacobian, 
the various pieces that cross don't really interact with each other.
In a soap film, two roughly perpendicular surfaces would merge
and probably create a singularity like the ones that are described 
below. See Example 3.a and Figures 2, 3, 5 below. 

If we replaced $A(f)$ with the surface measure
(or the Hausdorff measure) of $f(D)$, which may be smaller
if $f$ is not one-to-one, we would probably get a much better
description of soap films, but also a much harder existence 
theorem to prove. We shall return to similar ways of stating
Plateau problems in Section 4.e, 
when we discuss sliding minimal sets.

\ssi{\bf 2.b. Reifenberg's homology problem}
\ss

The second approach that we want to describe is due
to Reifenberg [R1] (1960). Let us state things for a 
$d$-dimensional surface (so you may take $d=2$ for simplicity).
Here we do not want to assume any a priori smoothness for the solution, 
it will just be a closed set $E$.
Because of this, we shall define the area of $E$ to be its
$d$-dimensional Hausdorff measure.
Recall that for a Borel-measurable set $E \i \R^n$, the $d$-dimensional 
Hausdorff measure of $E$ is
$$
\H^d(E) = \lim_{\delta \to 0^+} \H^d_\delta(E),
\leqno (2.3)
$$
where
$$
\H^d_\delta(E) = c_d \inf 
\Big\{ \sum_{j\in \Bbb N} \diam(D_j)^d \Big\},
\leqno (2.4)
$$
$c_d$ is a normalizing constant, and
the infimum is taken over all coverings of $E$ by a 
countable collection $\{ D_j \}$ of sets, with 
$\diam(D_j) \leq \delta$ for all $j$.
Let us choose $c_d$ so that $\H^d$ coincides with 
the Lebesgue measure on subsets of $\R^d$.

The main point of $\H^d$ is that it is a measure defined for
all Borel sets, and we don't lose anything anyway because $\H^d(E)$ 
coincides with the total surface measure of $E$ when
$E$ is a smooth $d$-dimensional submanifold. Also, we don't assume $E$ 
to be parameterized, which does not force us to worry about 
counting multiplicity when the parameterization is not one-to-one, 
and we allow $E$ to take all sorts of shapes, even if our boundary
is a nice curve; see Section 3 for examples of natural minimizers
that are not topological disks.

So we want to minimize $\H^d(E)$ among all closed sets $E$
that are bounded by a given $(d-1)$-dimensional set $\Gamma$.

\ms
We still need to say what we mean when we say that $E$
bounded by $\Gamma$, and Reifenberg proposes to define this
in terms of homology in $E$. He says that $E$
bounded by $\Gamma$ if  the following holds.
First, $\Gamma \i E$. This is not too shocking,
even though we shall see later soap films $E$ bounded by 
a smooth curve, but that seem to leave it at some singular point 
of $E$. See Figures 4, 5, and 17.
And even in this case, we are just saying that we see $\Gamma$
itself as (a lower dimensional) part of the film.
The main condition concerns the $\check{\rm C}$ech homology
of $E$ and $\Gamma$ on some commutative group $G$. Since
$\Gamma \i E$, the inclusion induces a natural homomorphism 
from $\check H_{d-1}(\Gamma;G)$ to $\check H_{d-1}(E;G)$, and
Reifenberg demands this homomorphism to be trivial.
Or, we could take a subgroup of $\check H_{d-1}(\Gamma;G)$
and just demand that each element of this subgroup
be mapped to zero in $\check H_{d-1}(E;G)$.

In the simple case when $d=2$ and $\Gamma$ is a simple
curve, $\check H_{d-1}(E;G)$ is generated by a loop
(run along $\Gamma$ once), and when we require this
loop to be a boundary in $E$, we are saying that there 
is a way to fill that loop (by a $2$-dimensional chain
whose support lies) in $E$.

\ss
Reifenberg proved the existence of minimizers in all dimensions, 
but only when $G = \Bbb Z_2$ or $G = \R /\Bbb Z$. 
This is a beautiful (although apparently very technical) 
proof by ``hands``'', where haircuts are performed on 
minimizing sequences to make them look nicer and allow 
limiting arguments. 
Unfortunately, difficulties with limits force him to use 
$\check{\rm C}$ech homology (instead of singular homology,
for instance) and compact groups.

Reifenberg also obtained some regularity results for the solutions;
see [Re1], [Re2]. 

Later, F. Almgren [Al3] 
proposed another proof of existence, which works for more general 
elliptic integrands and uses integral varifolds. But I personally find
the argument very sketchy and hard to read.

Recently, De Pauw [Dp] 
obtained the $2$-dimensional case when $\Gamma$ is a finite union 
of curves and now $G = \Bbb Z$, with a proof that uses currents; 
he also proved that in that case the infimum for this problem 
is the same number as for the size-minimizing currents of 
the next subsection. 
But even then it is not known whether one can build
size-minimizing currents supported on the sets that 
De Pauw gets. 

Reifenberg's solutions are nice and seem to give a good description
of many soap films. Using finite groups $G$ like $\Bbb Z_2$,
one can even get non-orientable sets $E$.
But there are some ``real-life" soap spanned by a curve 
that cannot be obtained as Reifenberg solutions.
See Example 3.c, for instance.

Anyway, many interesting problems (existence for other homologies 
and groups like $\Bbb Z$, equivalence with other problems)
remain unsolved in the Reifenberg framework of this subsection. 
The author does not know whether too much was done 
after [Re2], 
concerning the regularity of the Reifenberg minimizers.
But at least the regularity results proved for the Almgren minimal sets
(see Section 5) are also valid for the Reifenberg minimizers; we 
check this in Section 7 for the convenience of the reader.

\ssi{\bf 2.c. Integral currents}
\ss

Currents provide a very nice way to solve two problems at the same time. 
First, they will allow us to work with a much more general class of 
objects, possibly with better compactness properties. 
That is, suppose we want to work with surfaces in a certain class
$\cal S$, typically defined by some level of regularity, and we want
to prove the existence of some $S \in {\cal S}$ that minimizes
(some notion of) the area $A(S)$ under some boundary constraints.
Then let $\{ S_k \}$ be a minimizing sequence, which means that
each $S_k$ lies in ${\cal S}$ and satisfies the boundary constraints,
and that $A(S_k)$ tends to the infimum of the problem.
We would like to extract a subsequence of $\{ S_k \}$ that converges,
but typically the ${\cal S}$-norms of $S_k$ will tend to $+\infty$,
and we will not be able to produce a limit set $S\in {\cal S}$.
Of course even if $S$ exists, we shall not be finished, because we
also need to check that $S$ satisfies the boundary constraints, and
that $A(S) \leq \lim_{k \to +\infty} A(S_k)$, but this is a different
story. 

So, in the same spirit as for weak (or more recently viscosity) 
solutions to PDE, we want to define Plateau problems in a rough setting, 
and of course hope that as soon as we get a minimizer, we shall be able to 
prove that it is so regular that in fact it deserves to be called a 
surface.

The second positive point of using currents is that even with these
rough objects, we will be able to define a notion of boundary, 
inherited from differential geometry, and thus state a Plateau problem.
We need a few definitions.

A $d$-dimensional \underbar{current} is nothing but a continuous 
linear form on the space of smooth $d$-forms. 
This is thus the same as a $d$-vector valued distribution. 
In fact, most of the distributions that will be used here are 
($d$-vector valued) finite measures, which are thus not too wild.

There are two main examples of currents that we want to mention here.
The first one is the \underbar{current} $S$ \underbar{of integration on
a smooth, oriented surface} $\Sigma$ of dimension $d$,
which is simply defined by 
$$
\langle S,\omega\rangle = \int_\Sigma \omega
\ \hbox{ for every $d$-form $\omega$.}
\leqno (2.5)
$$
But of course we are interested in more general objects.
The second example is the
\underbar{rectifiable current} $T$ defined on a 
$d$-dimensional rectifiable set $E$ such that $\H^d(E) < +\infty$, 
on which we choose a measurable orientation $\tau$ and an 
integer-valued multiplicity $m$. 
Recall that a rectifiable set of dimension $d$ is a set $E$ such that
$\dsp E \i N \cup \bigcup_{j\in \Bbb N} G_j$, where $\H^d(N) = 0$
and each $G_j$ is a $C^1$ embedded submanifold of dimension $d$.
But we could have said that $G_j$ is the Lipschitz image 
of a subset of $\R^d$, and obtained an equivalent definition.
We shall only consider Borel sets $E$, and such that
$\H^d(E) < +\infty$. For such a set $E$ and $\H^d$-almost
every $x\in E$, $E$ has what is called an approximate tangent $d$-plane 
$x+V(x)$ at $x$, which of course coincides with the usual tangent plane
in the smooth case. A measurable orientation can be defined
as the choice of a simple $d$-vector $\tau(x)$ that spans $V(x)$, 
which is defined $\H^d$-almost everywhere on $E$ and measurable. We set
$$
\langle T,\omega \rangle = \int_E m(x) \; \omega(x)\cdot\tau(x) 
\, d\H^d(x)
\leqno (2.6)
$$
when $\omega$ is a (smooth) $d$-form, and 
where the number $\omega(x)\cdot\tau(x)$ is defined in a natural
way that won't be detailed here. We assume that the multiplicity
$m$ is integrable against ${\bf 1}_E d\H^d$, and then the integral
in (2.6) converges.

Return to general currents.
The \underbar{boundary} of any $d$-dimensional current $T$ 
is defined by duality with the exterior derivative $d$ on forms, 
by
$$
\langle \partial T,\omega \rangle = \langle T, d\omega \rangle
\ \hbox{ for every $(d-1)$-form $\omega$.}
\leqno (2.7)
$$

When $\Sigma$ is a smooth oriented surface
with boundary $\Gamma$, $S$ is the current of integration on 
$\Sigma$, and $G$ is the current of integration on $\Gamma$,
Green's theorem says that $\partial S = G$. 
This allows us to define Plateau boundary
conditions for currents. We start with a $(d-1)$-dimensional 
current $S$, for instance the current of integration on
a smooth oriented $(d-1)$-dimensional surface without boundary,
and simply look for currents $T$ such that
$$
\partial T = S.
\leqno (2.8)
$$
Notice that $\d\d = 0$ among currents, just because $dd=0$
among forms. So, if we want (2.8) to have solutions, we need
$\partial S = 0$; this is all right with loops (when $d=2$),
and this is why we want the smooth surface above to have no 
boundary. But other choices of currents $S$ would be possible.

We often prefer to restrict to integral currents.
An \underbar{integral current} is a rectifiable current $T$ 
as above (hence such that the multiplicity $m$ is integer-valued
and integrable against ${\bf 1}_E d\H^d$), and such that
$\d T$ is such a rectifiable current too. If $T$ solves (2.8),
the condition on $\d T$ will be automatically satisfied, 
just because we shall only consider rectifiable boundaries $S$.
The fact that we only allow integer multiplicities should make our
solutions more realistic, and make regularity theorem easier to prove.
Otherwise, we would expect to obtain solutions with very low density,
or (at best) foliated and obtained by integrating other solutions.

So here is how we want to define a Plateau problem in
the context of integral currents: we start from a given
integral current $S$, with $\d S= 0$ (for instance, the current of
integration on a smooth loop)
and we minimize the area of $T$ among integral currents $T$ such that
$\partial T = S$. 

\ms
The most widely used notion of area for currents is the 
\underbar{mass} $Mass(T)$, which is just the norm of $T$, 
seen as acting on the vector space of smooth forms 
equipped with the supremum norm (the norm of uniform convergence).
In the case of the rectifiable current of (2.6),
$$
Mass(T) = \int_E |m(x)| d\H^d(x).
\leqno (2.9)
$$
The corresponding Plateau problem works like a geometric measure 
theorist's dream. First, the problem of finding an integral current
$T$ such that $Mass(T)$ is minimal among all the solutions of (2.8)
has solutions in all dimensions, and as soon as $S$ is an integral
current with compact support and such that $\d S =  0$. 
This was proved a long time ago by Federer and Fleming [FF], [Fe1]; 
De Giorgi also had existence results in the codimension 1 case, 
in the framework of BV functions and Caccippoli sets.

What helps us here is the the lowersemicontinuity of $M$
(not so surprising, it is a norm), and the existence of a beautiful
compactness theorem that says that under reasonable circumstances 
(the masses of the currents $T_k$ and $\partial T_k$ stay bounded, 
their supports lie in a fixed compact set), the weak limit of 
a sequence $\{ T_k \}$ of integral currents is itself an integral current.

Moreover, the solutions of this Plateau problem (we shall call them
\underbar{mass minimizers}) are automatically very regular away from the boundary.
Let $T$ be a mass minimizer of dimension $d$ in $\R^n$, and denote by
$F$ its support (the closure of the subset of $E$ above where 
$m(x)$ is nonzero), and by $H$ the support of $S = \d T$. 
If $d= n-1$ and $n \leq 7$, $F$ is a $C^\infty$
embedded submanifold of $\R^n$ away from $H$,
and if $d=n-1$ but $n \leq 8$, $E$ may have a singularity set of dimension
$n-8$ away from $H$, but no more. 
See [Fe2]. 
And in larger codimensions, the dimension of the singularity set is
as most $d-2$ [Al5]. 
There were also important partial results by 
W. Fleming [Fl],  
Simons, and Almgren [Al2]; 
see [Mo5], Chapter 8 for details. 

With this amount of smoothness, mass minimizers cannot give a good
description of all soap films in $3$-space, because some of these 
have obvious interior singularities. Also, the fact that the notion of
current naturally comes with an orientation is a drag for some
examples (like M\"obius films). We shall discuss this a little 
more in the next section.

For the author, the main reason why mass minimizers do not seem 
to be a good model for soap films (regardless of their obvious 
mathematical interest), is because the mass is probably not 
the right notion of area for soap films. So we may want to consider 
the \underbar{size} $Size(T)$ which, is the case of the rectifiable current of 
(2.6), is defined by
$$
Size(T) = \H^d\big(\big\{ x\in E \, ; \, m(x) \neq 0 \big)\big\}.
\leqno (2.10)
$$
That is, we no longer count the multiplicity as in (2.9), we just
compute the Hausdorff measure of the Borel support. This setting
allows one to recover more examples of soap films, and to eliminate
some sets that are obviously not good soap films; see the next section.
But the price to pay is that the existence and regularity results
are much harder to get.

The Plateau problem that you get when you pick a rectifiable
current $S$ with $\d S = 0$, and try to minimize $Size(T)$
among all integral currents $T$ such that $\d T = S$, is very 
interesting, but as far as the author knows, far from being
solved. There are existence results, where one use intermediate
notions of area like $\int_E |m(x)|^\alpha d\H^d(x)$
(instead of (2.9) or (2.10)) for some $\alpha \in (0,1)$;
see [DpH]. 
But for instance we do not have an existence result when $d=2$
and $S$ is the current of integration over a (general) smooth closed
curved in $\R^3$.
The compactness theorem above does not help as much
here, because if $\{ T_k \}$ is a minimizing sequence, we control
$Size(T_k)$ but we don't know that $Mass(T_k)$ stays bounded,
so we may not even be able to define a limit which is a current.
See [Mo1], though, for the special case when $\Gamma$ is contained
in the boundary of a convex body.

The size-minimizing problem is not so different from the Reifenberg
Plateau problem of Section 2.b, and in some cases the infimum for
the two problems was even proved to be the same [Dp]. 

Let us not comment much about the (mostly interior) regularity results 
for size minimizers, and just observe that the regularity results for 
Almgren minimal sets apply to the support of $T$ when $T$ is a size 
minimizer. See Section 7.

\msi{\bf 3. SIMPLE CLASSICAL EXAMPLES}
\ss

It is probably time to rest a little, and try the various 
definitions above on a few simple examples. Most of the examples
below are very well known; hopefully they will be convincing, even 
though we can almost never prove that a given set, current, or surface,
is minimal. (Sadly, the main way to prove minimality is by exhibiting 
a calibration for the given object, which would often implies an algebraic
knowledge that don't have). We will try to give a more detailed 
account than usual of what happens, to make it easier for the reader.

\ssi{\bf 3.a. Crossing surfaces: the disk with a tongue}
\ss 
The next example is essentially the same as in [Al1].
Construct a (smooth) boundary curve $\Gamma$, 
as in Figure 1,  
which contains a main circular part, and two roughly parallel 
lines that cross the disk, near the center;
then solve a Plateau problem or drop the wire into a soap
solution and pull it back. 

Most probably (but in fact the author can't compute!), 
the parameterized solution of Douglas is a smooth, 
immersed surface $E_d$, that crosses itself along a 
curve $I$, which lies near the center of 
the disk and connects the two parallel lines;
the two pieces have no reason to really interact 
(the functional $A$ in (2.1) does not see that the two pieces
get close to each other). See Figure 1. 

Next orient $\Gamma$, let $S$ be the current of integration
on $\Gamma$, and minimize $Mass(T)$ among integral currents
$T$ such that $\d T = S$; then let $E_m$ denote the closed 
support of $T$. By Fleming's regularity theorem [Fl], 
$E_m$ is smooth away from $\Gamma$, so we can be sure that 
$E_m \neq E_d$ (also see Subsection 3.d for a simpler argument). 
Quite probably, 
$E_m$ looks like the surface suggested in Figure 2,  
which we could obtain from $E_d$ by splitting $E_d$ along $I$
into two surfaces, and letting them go away from each other
and evolve into something minimal. This creates a hole near
$I$, where one could pass without meeting $E_m$.
This also changes the topology of the surface: now $E_m$ is 
(away from $\Gamma$) a smooth surface with a hole, not the continuous
image of a disk. Also, $E_m$ is oriented,
because it is smooth and comes from the current $T$, whose
boundary lives in $\Gamma$. That is, locally and away from
$\Gamma$, we can write $T$ as in (2.6), where $\tau(x)$ is smooth;
then the fact that $\d T = 0$ locally implies that the multiplicity
$m$ is locally constant, and this gives an orientation on $E_m$
(if $E_m$ connected, as in the picture).

\vskip0.4cm
\centerline{ 
\includegraphics[height=50mm]{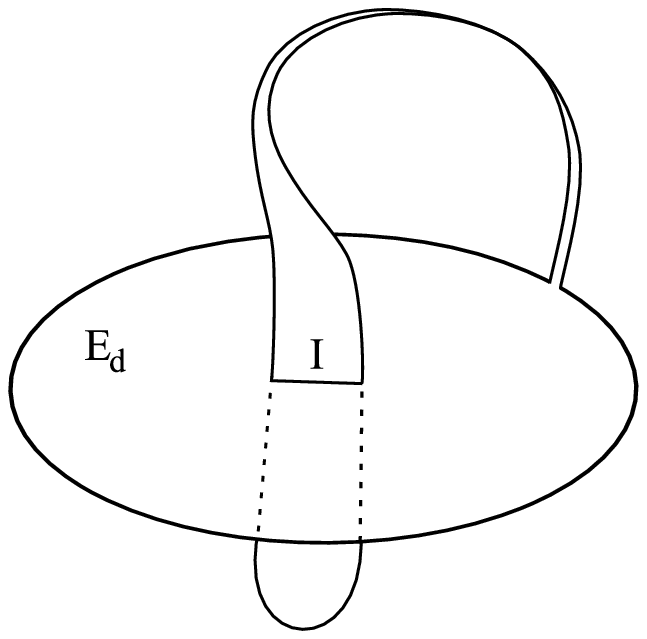}
\hskip1.5cm
\includegraphics[height=50mm]{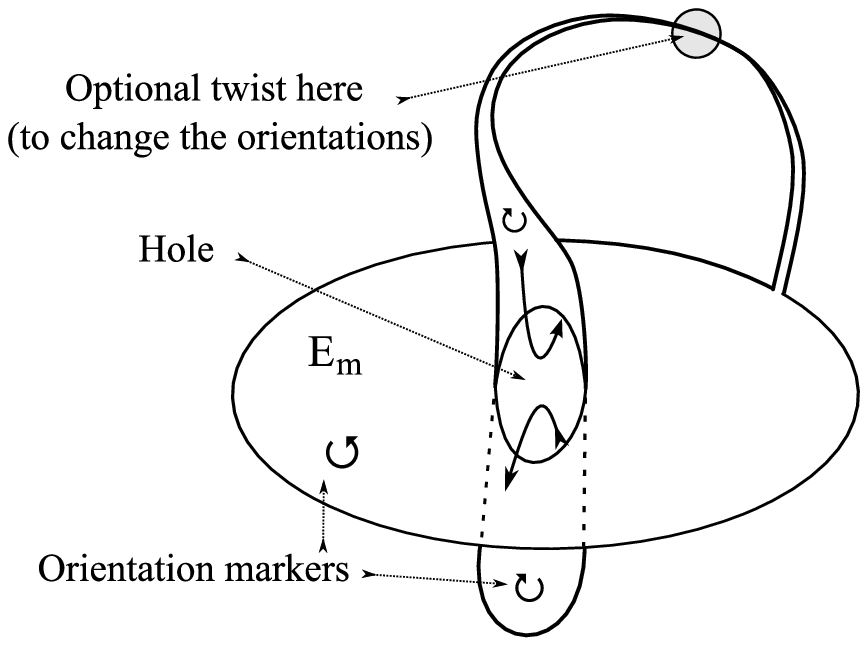}
}\medskip
\centerline{{\bf Figure 1.} The Douglas solution.
\hskip0.5cm {\bf Figure 2.} The mass minimizer $E_m$.}
\pari\noindent 
The circular little arrows mark the orientation, and the longer
oriented arrows try to show how $E_m$ turns in the two curved 
transition zones.
\medskip

Notice that the surface in Figure 2 
is orientable, but if we had chosen to use the symmetric way 
to split $E_d$ along $I$, we would have drawn a surface that goes in 
the direction of the front when we go down along the upper part 
of the tongue. 
This would have created a surface $\wt E_m$ which, near $I$, is roughly no one
the image of $E_m$ by the symmetry with respect to a vertical plane.
But $\wt E_m$ is not orientable, hence the mass minimizing Plateau
problem is not allowed to chose $\wt E_m$. The situation is reversed 
when we twist the thin part of $\Gamma$ to change the orientation, 
which physically should not matter much though.

Now let us build $\Gamma$ and plunge it into a film solution.
Based on a few rough experiments, the author claims the following
(which should surprise no one anyway).

\vskip.5cm 
\centerline{
\includegraphics[height=45mm]{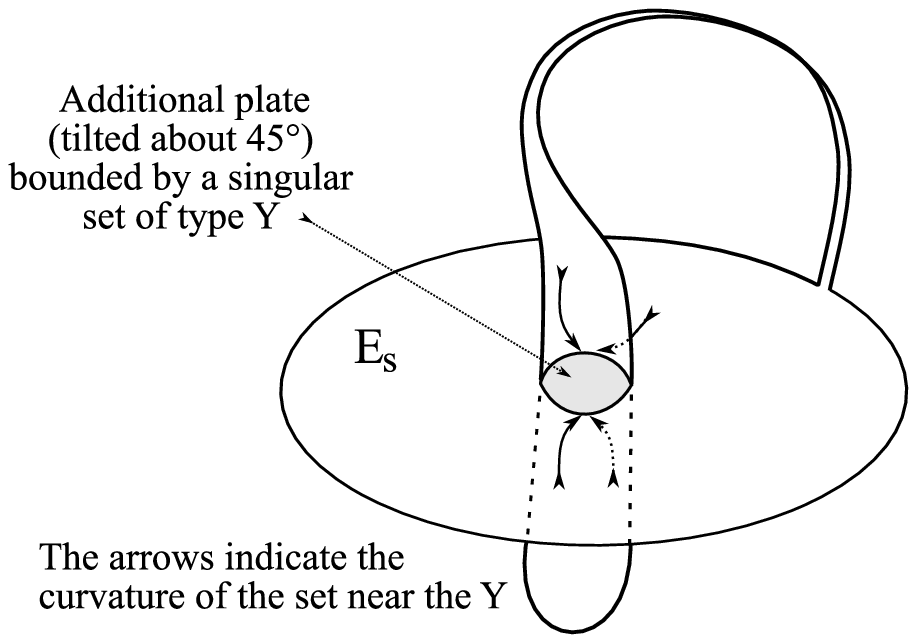}
\hskip0.2cm
\includegraphics[height=45mm]{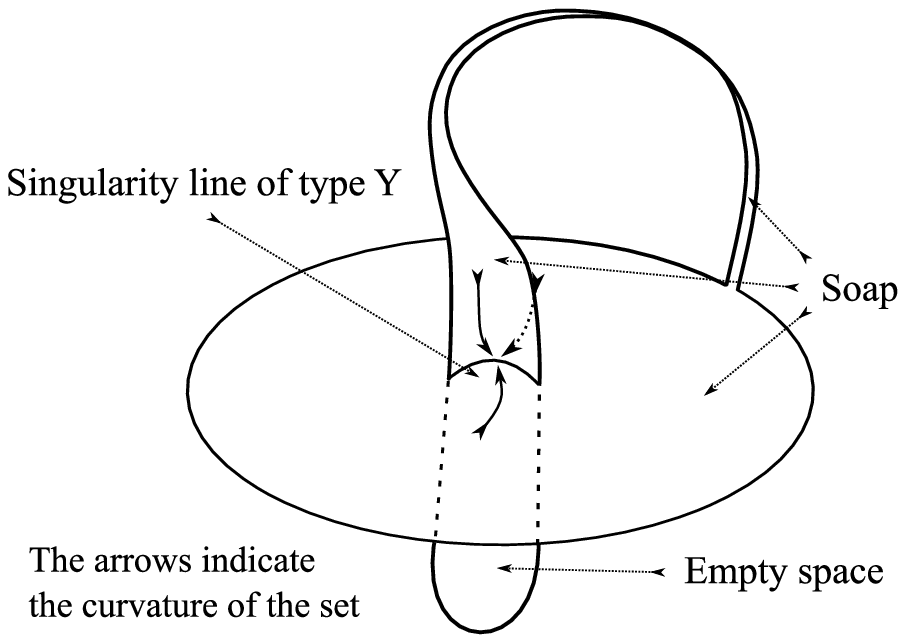} 
}\medskip
\centerline{{\bf Figure 3 (left).} The size minimizer $E_s$.}
\par
\centerline{{\bf Figure 4 (right).} 
A soap film which does not touch the whole boundary.}
\medskip 

We should not even try to obtain $E_d$, if we believe
that Almgren almost minimal sets give a good description of
soap films, and indeed the author managed not to see $E_d$.
Films that look like $E_m$ are not so hard to get; 
we also get a third set $E_s$ which we could roughly
obtain from $E_m$ as follows: add a small (topological) disk 
to fill the hole, and let again evolve into something minimal
(so that the boundary of the disk will be a set of singularities of type $Y$); 
see Figure 3, or the left part of Figure 5. 
Also see the double disk example in Figure 7 
for a more obvious singularity set of type $Y$).
Or we could obtain $E_s$ from $E_d$ by pinching it along $I$.
In practice, we often obtain $E_s$ first, especially if our
tongue is far from perpendicular to the horizontal disk
(otherwise, it seems a little less stable), and one can 
obtain $E_m$ from $E_s$ by killing the small connecting disk.

We can also kill the lower part of the tongue by touching it,
and get a singular curve of type $Y$ where the upper tongue connects
to the large disk-like piece. See Figure 4, or the right-hand   
part of Figure 5.
When this happens, we get an example of a soap film 
whose apparent boundary is a proper subset of
the curve $\Gamma$. This phenomenon is known, 
and has been described in the case of a trefoil knot
in [Br4], for instance. 
See Figure 17 
for another, more beautiful example.
For all these physical observations, the fact that the two main pieces 
of the boundary are connected and make a single curve is not useful, 
and we get a similar description when $\Gamma$ is composed of 
two disjoint non parallel ellipses (a large one and a small one) 
with the same center.

\centerline{
\includegraphics[height=50mm]{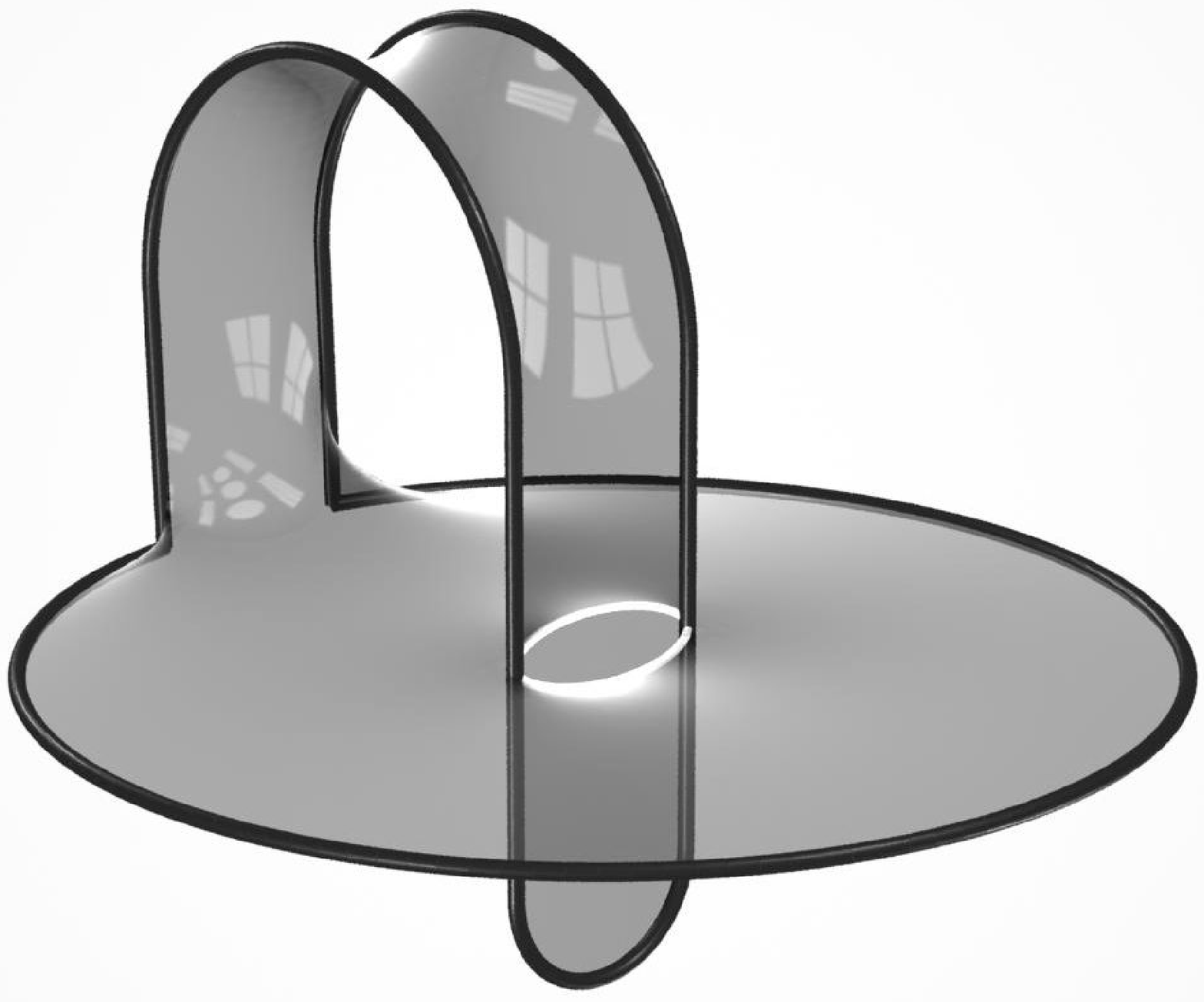}
\hskip0.8cm
\hskip 0.7cm
\includegraphics[height=50mm]{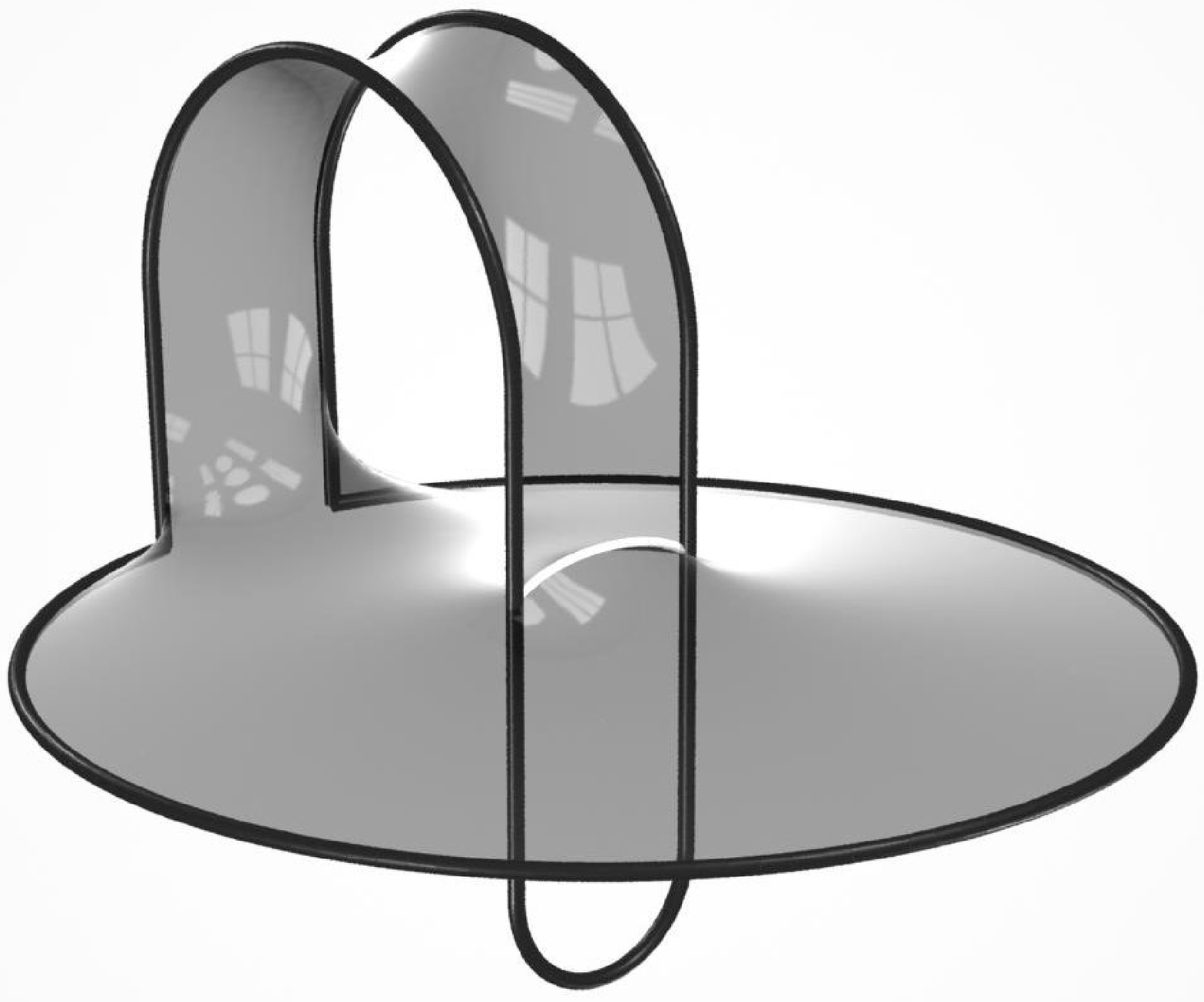} 
} \vskip-0.5cm
\noindent
{\bf Figure 5.} The sets of Figure 3 and 4. 
Images by John M. Sullivan, Technische Universit\"at Berlin,
used by permission.
\medskip

Finally (returning to the case of a single curve), 
if the part that connects the tongue to the 
circular part of $\Gamma$ is not too long or complicated, 
we obtain an often more stable, more complicated soap
film with a singularity of type $T$, like the set depicted
Figure 6. 
In fact, this one is often easier to get, and we
then retrieve $E_s$ and $E_m$ by removing faces.
See Section 3.e, and in particular Figures 9 and 10, for explanations 
and pictures of type $T$ singularities.

We leave it as an exercise for the reader to determine
whether $E_s$ (or its vaguely symmetric variant $\wt E_s$)
is probably the support of a size minimizer with the boundary $S$,
and whether $E_s$ and $\wt E_s$ are probable solutions
of Reifenberg's problem. We shall treat the simpler example of 
two circles instead.

Let us nonetheless say that (not even based on real experiments) 
the author believes that the fact that soap films will choose $E_m$ 
rather than $\wt E_m$, mostly depends on the angle of the tongue 
with the horizontal disk, or rather the way we pull $\Gamma$
out of the soap solution, and is not based on orientation. 
Often we get $E_m$ or $\wt E_m$ by passing through $E_s$ and 
$\wt E_s$ first, which themselves should not depend much on orientation.

\vskip0.5cm 
\centerline{
\includegraphics[height=50mm]{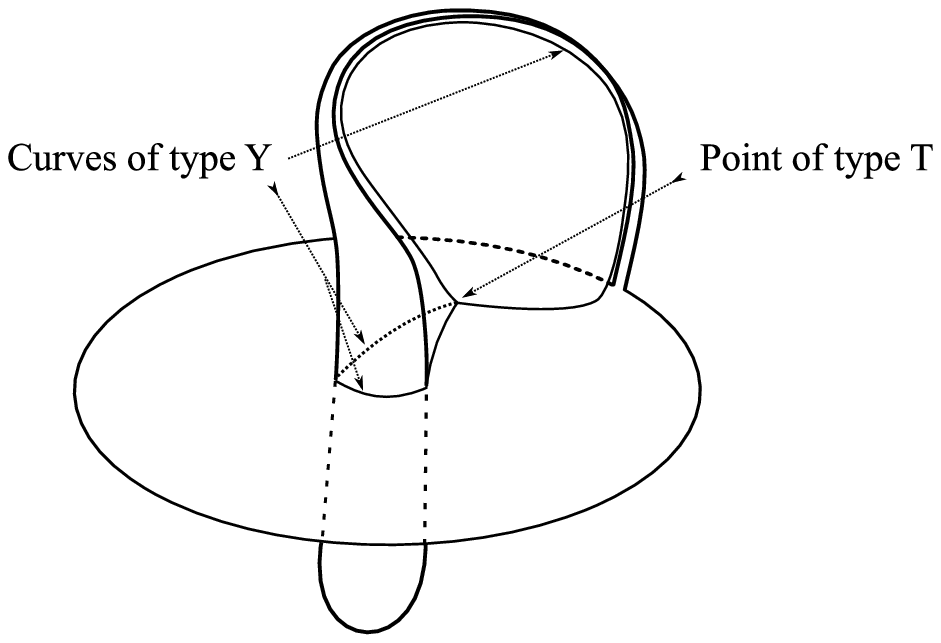}
}
\medskip\noindent
{\bf Figure 6.} Another soap film 
bounded by the same curve. The picture only shows the singular set 
composed of curves where $E$ looks like a $Y$, 
and which meet at a point where $E$ looks like a $T$.
\medskip

\ssi{\bf 3.b. Two parallel circles}
\ss
Let $\Gamma$ be the union of two parallel circles $C_1$ and $C_2$,
as in Figure 7. 

There are three obvious soap film solutions: a catenoid $H$
(that here looks a lot like a piece of vertical cylinder),
the union $D_1 \cup D_2$ of two disks, and a set $E$
composed of a slightly smaller disk $D$ in the center,
connected along the circle $\d D$ to $C_1$ and $C_2$ by two piece of 
catenoids that make angles of $120$ degrees along $\d D$.
See Figure 7. 

For the analogue of the Douglas problem, we would probably decide
first to parameterize with two disks or with the cylinder, and then
get $D_1 \cup D_2$ or $H$.

For the mass minimizing problem, we first choose an orientation on 
both circles $C_1$ and $C_2$, and for instance,
use the sum $S = S_1 + S_2$, where $S_i$ is the current of integration
on $C_i$.
If we choose parallel orientations, $H$ is not allowed 
because the orientations do not fit, and, with a fairly easy 
proof by calibration (use the definition of $\d T$
to compute $\langle T, dx_1 dx_2 \rangle$)
we can show that the current of integration on $D_1 \cup D_2$ 
is the unique mass minimizer.

If we choose opposite orientations on the $C_i$ 
(or equivalently set $S = S_1 - S_2$), then both $D_1 \cup D_2$ and $H$ 
will provide acceptable competitors,
and we should pick the one with the smallest mass. 
Here too, we should be able to pick a vector field
(pick one for each of the two competitors) and use it to
prove minimality by a calibration argument.

We could also try another choice of $S$,
like $S= S_1 +  2S_2$, but we would get the set $D_1 \cup D_2$
again.

The size-minimizing problem has a different solution
when we take parallel orientations. Let $D$ be a 
slightly smaller disk, parallel to $D_1$ and $D_2$,
and that we put right in the middle. Complete $D$ with two
pieces of catenoids $H_1$ and $H_2$, connecting $\d D$ 
to the $C_i$ (see Figure 7). 
The best choice will be when the $H_i$ make an angle
of $120$ degrees with $D$ along $\d D$.

It is easy to construct a current $T$, with
$\d T = S$, and which is supported on $E$: take twice
the integration on $D$, plus once the integration on
each catenoid, all oriented the same way.
The fact that $\d T = S$ is even easier to see when
one notices that $T = T_1 + T_2$, where $T_i$ is
the current of integration on $D \cup H_i$ and 
$\d T_i = S_i$. If the two circles are close enough,
it is easy to see that $Size(T)$ is smaller
(almost twice smaller) than the size of the mass 
minimizer above. It should not be too hard to show that
in this case, $T$ is the unique size minimizer,
but the author did not try.
Again, $H$ is not allowed, because the orientations 
do not fit.

\vskip.5cm 
\centerline{
\includegraphics[height=28mm]{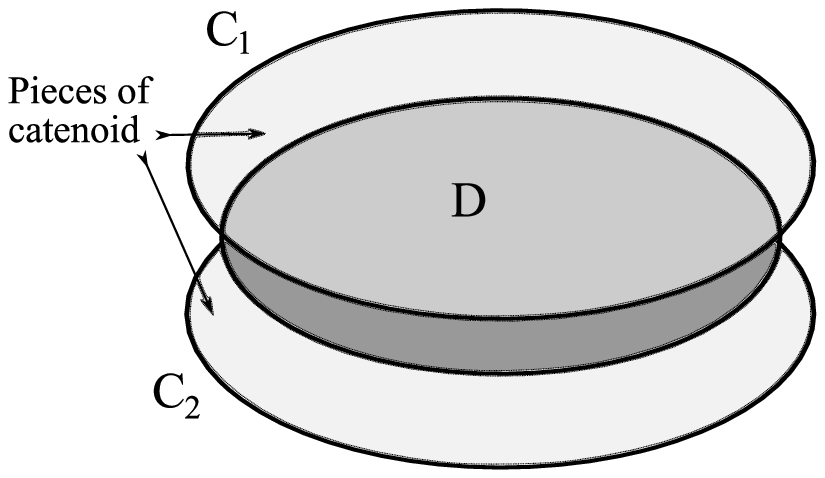}
\hskip1.2cm
\includegraphics[height=30mm]{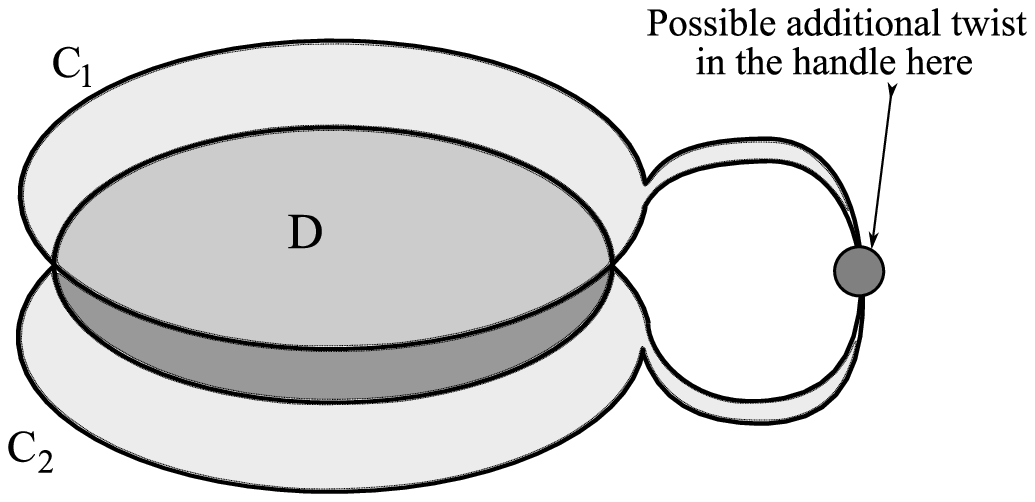} 
}\medskip
\centerline{{\bf Figure 7.} The set $E$. \hskip0.9cm
{\bf Figure 8.}  Same thing with a handle.}
\medskip 

When we use opposite orientations on the $C_i$,
$H$ is allowed, and $E$ never shows up: if $T$
is supported on $E$ and $\d T = S_1-S_2$,
the multiplicity on each $H_i$ should be constant,
equal to $(-1)^i$ (check $\d T$); the multiplicity
on $D$ should also be constant, and in fact equal
to $0$ if we want $\d T$ to vanish near $\d D$.
Then the support of $T$ is $H_1 \cup H_2$, and in
fact $H$ was doing even better than both 
$H_1 \cup H_2$ and $E$, so we should have no regret.

Returning to the case of parallel orientations, 
notice that when $H^2(E) \leq \H^2(D_1 \cup D_2)$,
$T$ is a size minimizer but $Size(T) = \H^2(E) > \H^2(H)$,
which is possible because $H$ is not a competitor.

The situation is simple enough for us to comment on the
Reifenberg problem. But even then, let us just think
about simplicial homology on the group $\Bbb Z$.
Here $\Gamma = C_1 \cup C_2$, and there are two obvious
generators $\gamma_1$ and $\gamma_2$ for the homology group 
$H_1(\Gamma, {\Bbb Z})$. Let us choose the $\gamma_i$ so
that they correspond to parallel orientations of the $C_i$.

One way to state a Plateau problem is to minimize $\H^2(F)$
among sets $F$ such that
$\gamma_1+\gamma_2$ represents a null element
in the homology group $H_2(F,{\Bbb Z})$. 
Notice then that $F = D_1 \cup D_2$ is allowed
(because $D_1 \cup D_2$ gives a simplicial chain
supported in $F$ and whose boundary is $\gamma_1+\gamma_2$),
but $H$ is not. To check this we would have to check that
$\gamma_1 + \gamma_2$ does not vanish in $H_2(H,{\Bbb Z})$,
which we kindly leave as an exercise because we do not want 
to offend any reader that would know anything about homology.
Here $E = D \cup H_1 \cup H_2$ is allowed too, 
because each $D \cup H_i$ contains a simplicial chain whose
boundary is $\gamma_i$. And, if $D_1$ and $D_2$ are close
enough, the Reifenberg minimizer will be $E$ (but the verification 
will be more painful than before). 

The situation would be the same if we required that
$\gamma_1+2\gamma_2$ represents a null element
in $H_2(F,{\Bbb Z})$, or if we required the whole group
$H_1(\Gamma,{\Bbb Z})$ to be mapped to $0$ in 
$H_2(F,{\Bbb Z})$. 

On the other hand, if we just require $\gamma_1 - \gamma_2$
to be mapped to $0$, then both $H$, $D_1\cup D_1$, and
$E$ are allowed, and the minimizer will be $H$ if $D_1$ is close
to $D_2$. 

As far as the author knows, the verifications would be more painful, 
but at the end similar to the verifications for size-minimizing 
currents. The point should be that if $F$ is a competitor in the
Reifenberg problem, we should be able to construct chains inside
$F$ (whose boundary is a representative for $\gamma_1+\gamma_2$, 
for instance), possibly approximate them with polyhedral chains,
and integrate a differential form on them to complete a calibration 
argument. This is an easy exercise, but slightly above the 
author's competence.

We may also want to work on the group ${\Bbb Z}/2{\Bbb Z}$; 
then $\gamma_1+\gamma_2 = \gamma_1 - \gamma_2$, $H$ is allowed
in all cases, and is the unique minimizer when the disks are close
to each other.

Very easy soap experiments show that the three sets $E$,
$D_1 \cup D_2$, and $H$ are constructible soap films,
with a noticeable preference for $E$ when $D_1$ and
$D_2$ are close to each other.

\ssi{\bf 3.c. Two disks but a single curve}
\ss
In the previous example, we managed to represent all
the soap films as reasonable competitors in one of the 
two standard problems about currents, even though in some
cases, stable soap films are not absolute minimizers.
But we expect difficulties in general, again because of
orientation issues

In fact, is easy to produce a M\"obius soap film, 
which is not orientable and hence is not the support of a 
mass- or size-minimizing current. Unless we use some tricks
(like work modulo $2$, or use covering spaces), as in 
Sections 4.a and 4.b. For instance, the author believes
that a M\"obius soap film may be a solution of the Reifenberg 
problem above, where we decide to compute homology 
over the group $\Bbb Z_2$

Let us sketch another example, which is just a minor 
variant of the previous one. 
Cut two very small arcs out of $C_1 \cup C_2$, one above
the other, and replace them with two parallel curves $g_1$
and $g_2$ that go from $C_1$ to $C_2$; 
this gives a single simple curve $\Gamma_1$, 
as in Figure~8. 
This is also almost the same curve that is represented in K. Brakke's
site under the name ``double catenoid'' soap films.
There are many constructible soap films bounded by $\Gamma_1$
(see Brakke's site) but let us concentrate on the one that looks like 
the set $E$ above, plus a thin surface bounded by $g_1$ and $g_2$
(call this set $E_1$).
The difference with the previous example is that now there is only
one curve, the orientations on the two circular main parts
of $\Gamma_1$ are opposite, and (as in the case of opposite 
orientations above) $E_1$ is not the support of an integral
current $T$ such that $\d T$ is the current of integration on 
$\Gamma_1$. Similarly, $E_1$ cannot solve the Reifenberg problem
either, just because a strictly smaller subset is just as good
(we can check that we can remove the central disk from $E_1$,
and that $\Gamma_1$ is the boundary of the remaining surface, 
which is orientable).

If instead of taking $g_1$ and $g_2$ parallel, we make them
twist and exchange their upper extremities, we get parallel
orientations again, and the corresponding set $\wt E_1$
is a competitor in the Reifenberg and size-minimizing problems.

Two last comments for this type of examples: if we have only one Plateau
problem, as is the case in any of the categories above if the boundary is
just one curve, then in generic situations we can only get one solution. 
This is usually not enough to cover the variety of different soap films. 
But this alone would not be so bad if we could cover all the example as stable 
local minima. Also, orientation seems to be the main source of trouble here,
which suggests that we use varifolds. But boundary conditions are harder to 
define for varifolds.

\ssi{\bf 3.d. The minimal cones $Y$ and $T$, and why mass minimizers do not cross}
\ss

There are three types of tangent cones that can easily be seen
in soap films: planes (that we see at all the points where the
film is smooth), the sets $Y$ composed of three half planes 
that meet along a line with $120$ degrees angles, and the sets $T$ 
that are obtained as the cone over the union of the edges of a regular 
tetrahedron (seen from its center).
See Figure 9.  
Also see Figure 10 
for a soap film 
with a singularity of type $T$.

These are the three cones of dimension $2$ in $\R^3$ that are 
minimal sets in the sense of Almgren 
(see Section 5 below), but let us consider currents for the moment.

The plane is the only cone that can be seen as a blow-up limit
of the support of a $2$-dimensional mass minimizer in $\R^3$; 
this comes from Fleming's regularity theorem, but we want to 
discuss this a little more.

\vskip.5cm 
\centerline{
\includegraphics[height=36mm]{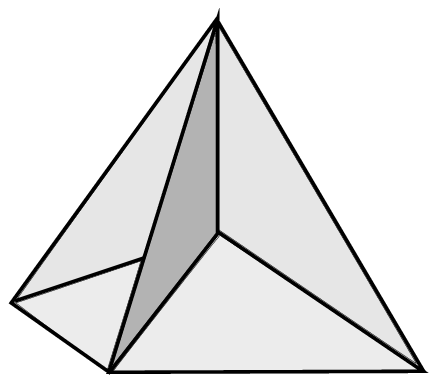} 
\hskip2.1cm
\includegraphics[height=34mm]{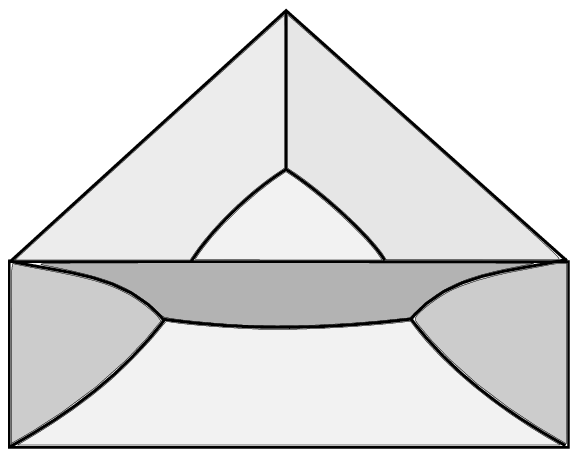} 
}
\medskip
\centerline{$\,$\hskip1.2cm
{\bf Figure 9.} The set $T$. \hskip1.2cm
{\bf Figure 10.}  A film with three $T$ points.}
\medskip

Let us first say why the support of a mass minimizer $T$ never 
looks like two smooth surfaces that cross neatly.
Denote by $E$ this support, and suppose that in a small box, 
$E$ looks like the union of the two planes (or smooth surfaces)
suggested in Figure 11. 
First assume that $T$ has multiplicity $1$ on these planes, 
with the orientation suggested by the arrows. Replace $E$,
near the center, with two smoother surfaces (as suggested in
Figure~12). 
This gives a new current $T_1$, and it is easy to see that
$\d T_1 = \d T$ (both $\d T$ and $\d T_1$ vanish near the
center, and the contributions away from the center are the same too).
It is also easy to choose the two smoother surfaces so that 
$Mass(T_1) < Mass(T)$.

\vskip.5cm 
\centerline{
\includegraphics[height=25mm]{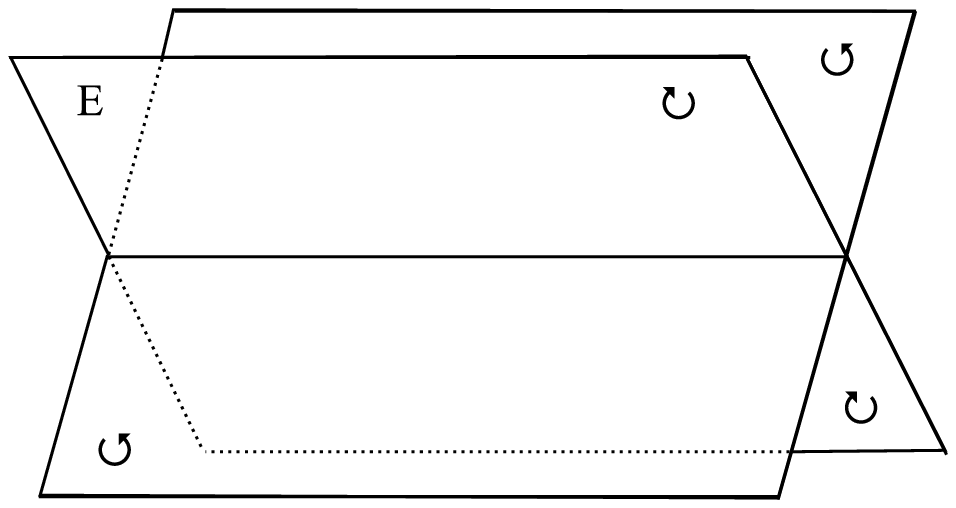} 
}
\smallskip
\centerline{
{\bf Figure 11.} Two oriented planes} 
\medskip

Notice that the modification suggested in Figure 13) 
is not allowed (how would we orient the surfaces?), and that the
modifications above are the same as what we suggested in
the tongue example of Subsection 3.a.

\vskip.5cm 
\centerline{
\includegraphics[height=33mm]{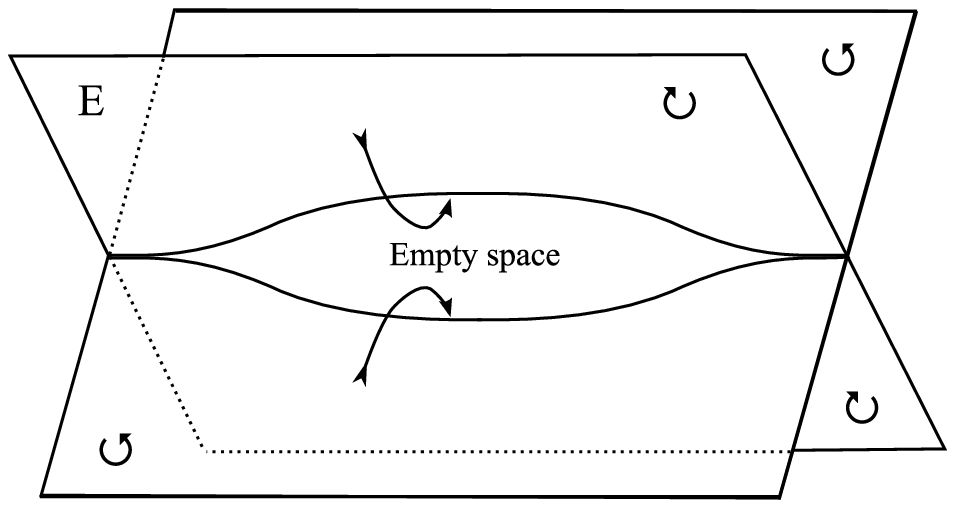}
\hskip1.1cm
\includegraphics[height=33mm]{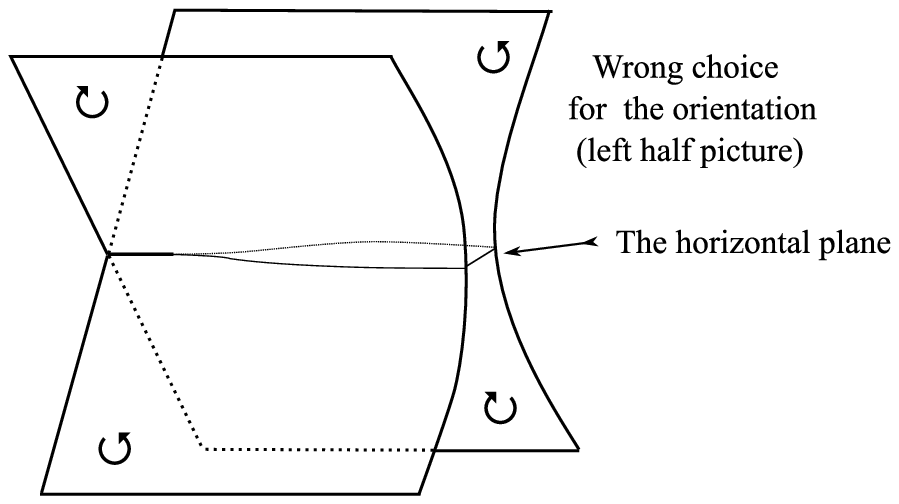}  
}
\medskip
\centerline{
{\bf Figure 12.} A better competitor.
\hskip0.2cm 
{\bf Figure 13.} Unauthorized competitor (half picture).
}
\medskip
We assumed above that $T$ has multiplicity $1$ on the two
planes. If the multiplicities are different, say, 
$1 \leq m_1 < m_2$, denote by $T_0$ the current on $E$
with the constant multiplicity $m_1$, replace $T_0$
with $m_1$ times the current of integration on the two
smooth surfaces, and keep $T-T_0$ as it is. We still get
a better competitor. This looks strange, and the author sees
this as a hint that the mass is too linear a functional to
be completely honest.

\ms
Return to the cone $Y$, and let us first say why it is the support 
of a current with no boundary. Put an orientation on the three 
faces $F_j$ that compose $Y$. 
The orientation of each $F_j$ gives an orientation
of the common boundary $L$, and we can safely assume that
the three orientations of $L$ that we get coincide.
Denote by $T_j$ the current of integration on $F_j$, and set
$T = T_1+T_2-2T_3$. 
It is easy to see that $\d T = 0$ because the three
contributions cancel.

One can show (and the best argument uses a calibration; 
see [Mo1]) that $T$  
is locally size minimizing. More precisely, for each choice of
$R > 0$, set $T_R = {\bf 1}_{B(0,R)} T$; the computation
shows that $\d T_R$ is of the form $S_1 + S_2 -2 S_3$,
where $S_j$ denotes the current of integration on the
(correctly oriented) half circle $F_j \cap \d B(0,R)$,
and one can show that for each $R > 0$,
$T_R$ is the unique size minimizer $W$ under the 
Plateau condition $\d W = \d T_R$. 

Now $Y$ is not the support of a local mass minimizer $T$:
the multiplicity would need to be constant on the three faces $F_j$,
the sum of these multiplicities should be zero because $\d T = 0$ 
near the line, and then we could split $T$ into two pieces
(each with two faces) that we could improve independently, 
as in the previous case.

\ms
The third minimal cone $T$ is also the support of a 
current $W$ with vanishing boundary; again we have to
put (nonzero integer) multiplicities $m_j$ on the six faces $F_j$,
so that their contribution to each of the four edges cancel.
Let us give an example of multiplicities that work.
Let $H$ denote the closed convex hull of the tetrahedron
that was used to define $T$. Then $T \cap \d H$ is
composed of six edges  $\Gamma_j = F_j \cap  \d H$;
we orient $\Gamma_j$ so that $\d ({\bf 1}_{H} T_j) = S_j$
along $\Gamma_j$, and then we set $S = \sum_j m_j S_j$.
Figure 14 
describes a choice of multiplicities and orientations
of the $\Gamma_j$ for which $\d S = O$ and so (this is the same
condition) $\d W = 0$.

With such multiplicities, ${\bf 1}_{H} W$ is the only size minimizer
for the Plateau condition $\d T = S$. See [Mo1]. 
Of course it is not a mass minimizer, because it contains a lot
of singularities of type $Y$, which are not allowed.

\vskip.5cm 
\centerline{
\includegraphics[height=35mm]{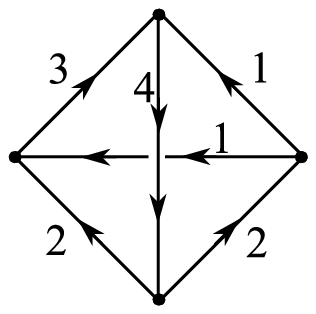}
\hskip1.0cm
\includegraphics[height=30mm]{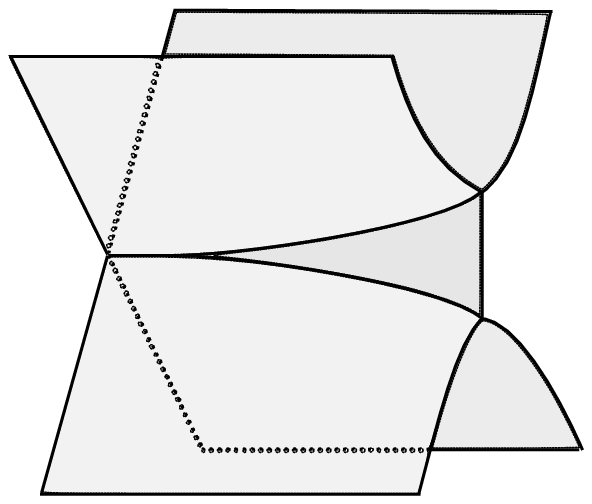} 
\hskip0.3cm
\includegraphics[height=30mm]{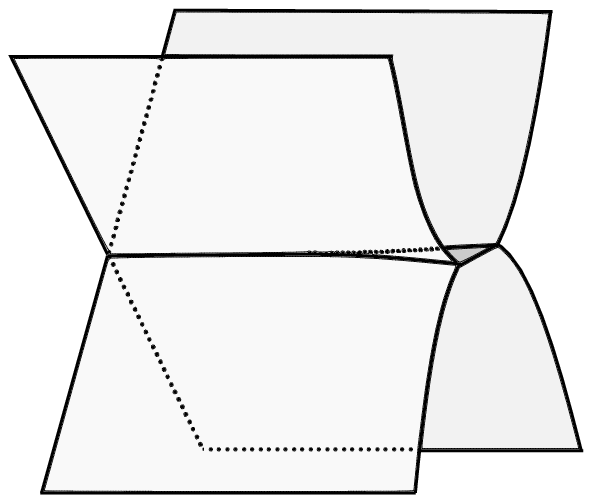}
}
\medskip
\noindent
{\bf Figure 14 (left).} Orientations and multiplicities
to make a size minimizer supported on $T$. The dots represent 
vertices of a tetrahedron, and the arrows sit on the edges $\Gamma_j$.
\par\noindent
{\bf Figure 15 (right).} Two size-competitors for the union
of two planes (half pictures). 
\medskip

Return to the union of two planes in Figure 11 
and let $T$ be obtained by putting nonzero integer multiplicities
on the four half faces; suppose that $\d T = 0$ in some large ball and, 
for convenience, that the two planes are nearly orthogonal. 
Then the two sets suggested in Figure~15  
support currents that have a smaller size than $T$. 
(It could be that we can take a vanishing multiplicity
on the middle triangular surface, and then one of the competitors of 
Figures~12 and 13 
is allowed and does even better.) 

\ssi{\bf 3.e. Mere local minima can make soap films}
\ss

We have seen that (real-life) soap films can exist, even if they do not
minimize mass, size, or Hausdorff measure in one of the Plateau
problems quoted above. If $\Gamma$ is a curve, each Plateau problem
usually comes with one solution, and many soap films exist, often
with different topologies (but maybe the same homology constraints).
We shall see in the next section a few attempts to multiply the 
number of Plateau problems, in order to accommodate more examples.

It is expected that soap films are not necessarily global 
minima, even for a given topology, i.e., that stable local
minima work as well. Things like this are unavoidable.
But the situation is even worse:
Figure 16  
shows a soap film which can be retracted, 
inside itself, into its boundary which is a curve.
But the retraction is long, and the soap will not see it 
and stay at the local minimum. The author admits he does not
see the retraction either, and was unable to construct
a soap example.

This example is due to J. Frank Adams (in the appendix of [Re1]), 
and the picture comes from K. Brakke's web page.

We add in Figure 17 
one of the many pictures of soap films bounded by a union
of three circles in a Borromean position. This one leaves one
of the three circles for some time, like the film of 
Figures 4 and 5. 
See the web page of K. Brakke for many more.

\vskip.1cm 
\centerline{
\includegraphics[height=50mm]{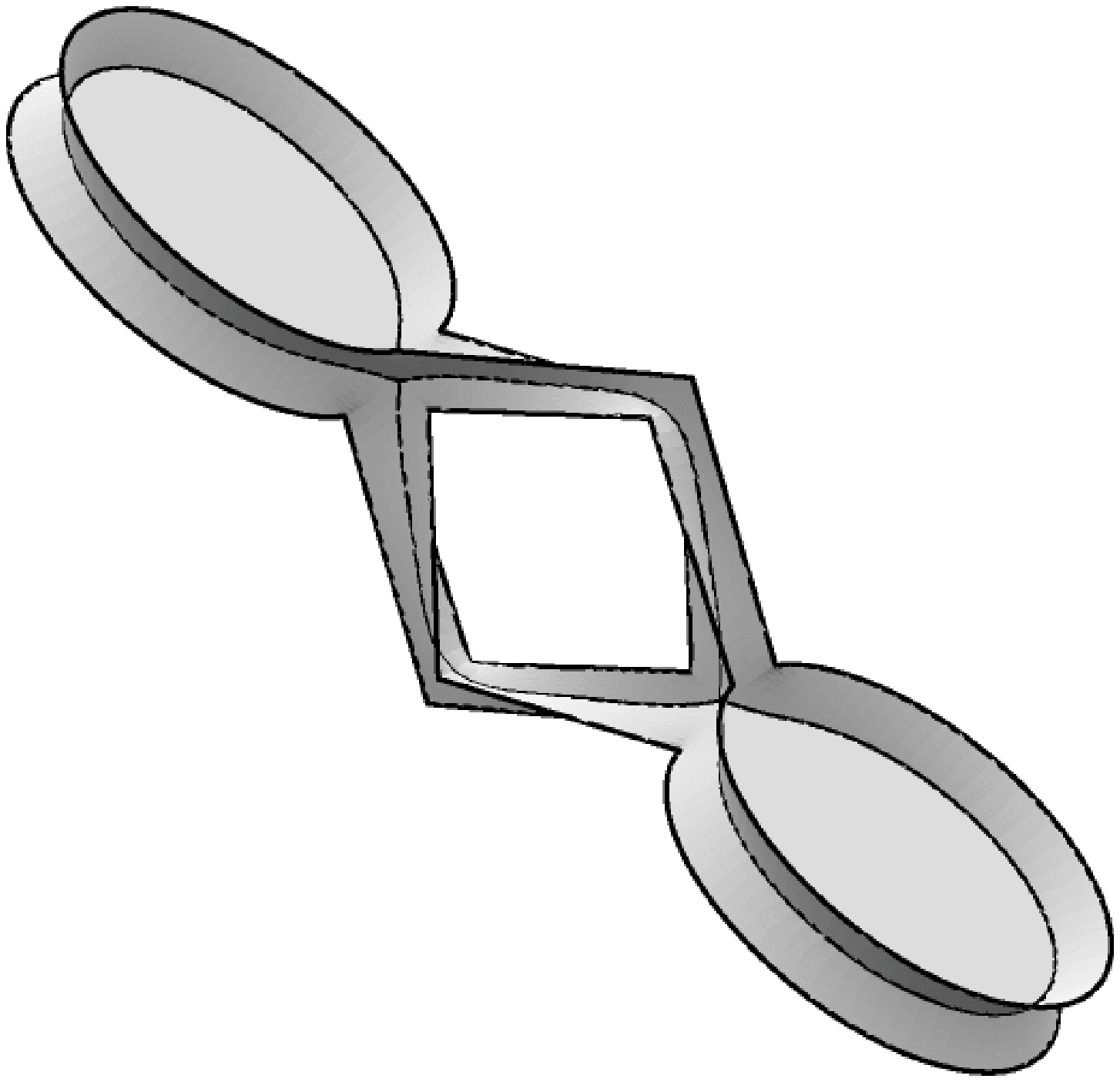}
\hskip1.0cm
\includegraphics[height=50mm]{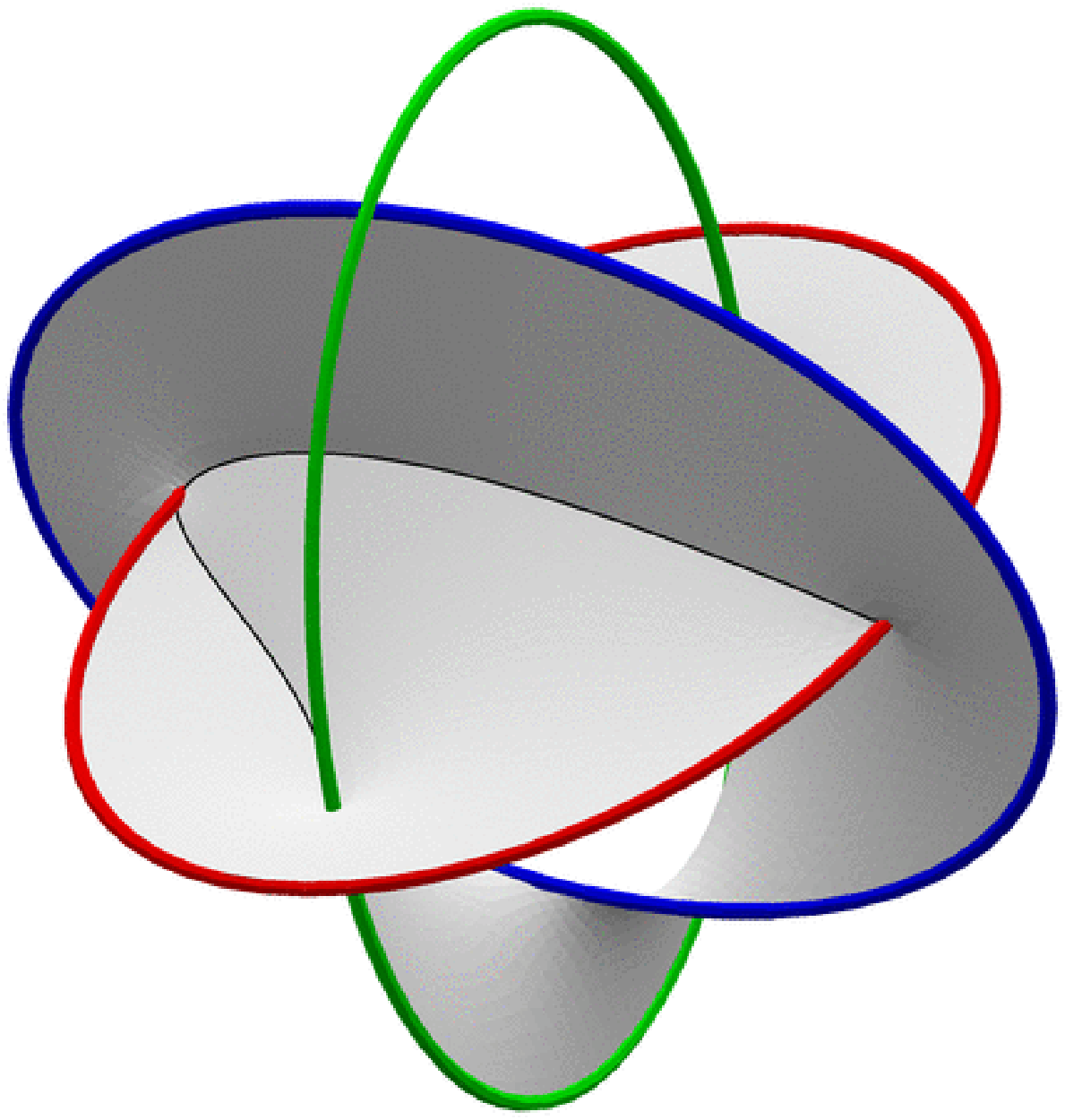}
}
\medskip
\noindent
{\bf Figure 16 (left).} A soap film that retracts onto itself
(J. Frank Adams).
\par\noindent
{\bf Figure 17 (right).} One of the many soap films
bounded by a Borromean ring.
Both pictures are courtesy of K. Brakke.
\medskip

\msi{\bf 4. OTHER MODELS FOR PLATEAU PROBLEMS}
\ss

Various tricks, often clever, have been invented to increase
the number of Plateau problems, for instance associated to a given curve, 
and thus accommodate the various examples. As we have seen, one of the 
unpleasant things that we often have to deal with is orientation.
We briefly report some of these tricks.

\ssi{\bf 4.a. Compute modulo $p$}
\ss
One can define integral currents modulo $k$, by saying that
two integral currents are equivalent when their difference is 
$k$ times an integral current. Then M\"obius strips, for instance, 
can can support currents modulo 2 without being orientable, and even
solve a mass minimizing Plateau problem as in Section 2, but with
integral currents modulo $2$. Similarly, the double disk set $E$
of Figure 2 is the support of a a current 
modulo 3 with multiplicity $\pm 1$.

\ssi{\bf 4.b. Use covering spaces}
\ss
In [Br4], K. Brakke manages to treat many 
of the simple soap film examples above in the general formalism of
mass (not even size!) minimizers. His construction works for 
sets of codimension $1$.
He starts from the base manifold $M = \R^n \sm \Gamma$, where
$\Gamma$ is our boundary set, constructs a covering space over
$M$, which branches along $\Gamma$, and eventually gets the
desired soap film $E$ as the projection of boundaries of
domains in the covering space, that minimizes the mass under
some constraints that we don't want to describe here. 
The boundaries are orientable even if $E$ is not,
and almost every point of $E$ comes from two boundary points 
in the covering (as if, a little as what happens in real soap film, 
$E$ were locally composed of two layers, coming from different
levels in the covering).

The construction is very beautiful, but apparently limited
to codimension 1. Also, the fact that a whole new construction 
seems to be needed for each example is a little unsettling.

\vfe 
\ssi{\bf 4.c. Use varifolds}
\ss

Almgren [Al3] 
used varifolds to present another proof of existence for
the solutions of the Reifenberg problem of Section 2.b.
Varifolds look like a very nice concept here, because they 
don't need to be orientable (which was an unpleasant
feature of currents), and we can compute variations of the
area functional on them. This gives a notion of stationary 
integral varifold, which englobes the minimal cones above,
and Allard and Almgren prove compactness theorems on
these classes that are almost as pleasant as for the
integral currents above. See [Al1], 
where the author expresses very high hopes that varifolds 
are the ultimate tool for the study of minimal sets and
Plateau problems. But (possibly because varifolds are not 
oriented), it seems hard to state and solve a Plateau problem, 
as we did in Section 2.c for currents with the boundary operator $\d$.
We shall return briefly to this issue in Section 6.

\ssi{\bf 4.d. Differential chains}
\ss

J. Harrison [Ha3] 
proposes yet another way to model and solve a Plateau problem.
Again the point is to get rid of the difficulties with orientation,
and at the same time to keep a boundary operator.
To make things easier, let us consider $2$-dimensional sets in $\R^3$. 
We are looking for a representation of soap films as slight generalizations of 
``dipole surfaces", where a dipole surface associated to a smooth
oriented surface $\Sigma$ is defined as 
$$
\Delta S = \lim_{t \to 0} {1 \over t} (S^+-S),
\leqno (4.1)
$$
where $S$ is the current of integration on $\Sigma$ and 
$S^+$ is the current of integration on the surface $\Sigma^+$ 
obtained from $\Sigma$ by a translation of $t$ in the 
normal direction. Thus dipole surfaces are currents of 
dimension $2$, but they act on forms by taking an 
additional normal derivative
(which means that they have one less degree of smoothness).

One of the ideas is that the orientation disappears
(the dipole  surface obtained from $\Sigma$ with the opposite
orientation is the same), which allows one to represent the 
branching examples of Section 3 as (limits of) dipole surfaces. 
The boundary $\d (\Delta S)$ is a dipole version
of $\d S$, taken in the direction of the unit normal to $\Sigma$
(along $\d \Sigma$).

J. Harrison works with the closure $\cal F$ of all finite sums
of $\Delta S$ as above,
with the norm coming from the duality with Lipschitz $2$-forms
(with the $L^\infty$ norm on the form and its derivative), 
starts from a nice curve $\Gamma$, chooses a smooth normal vector
field on $\Gamma$, uses it to define a dipole curve $G$ based as above
on the current of integration on $\Gamma$, and finally looks 
for a current $T \in \cal F$ that solves $\d T = G$, 
and for which the correct analogue of the mass 
(called the volume form) is minimal. 
In the case of a small $\Delta S$ above, the volume form
is equal to $\H^2(\Sigma)$; for general elements of $\cal F$, it is
defined by density.

It turns out that when we restrict to dipoles
(that is, limits of dipole polyhedral chains, for the norm 
of duality with the forms with one Lipschitz derivative),
there is a construction (based on filling the holes with a 
Poincar\'{e} lemma) that allows one to invert $\d$ 
and construct competitors in the problem above. 
This construction may be a main difference with the
situation of varifolds. 

But at the same time the elements of $\cal F$ are not so smooth, 
and one still needs to check some details, not only about 
the existence of minimizer for the problem above, but also 
to get some control on the solutions and show that they are
more than extremely weak solutions. For instance, 
it would be nice to prove that their support is a rectifiable set,
maybe locally Almgren-minimal, as in the discussion below.

\ssi{\bf 4.e. Sliding deformations and sliding Almgren minimal sets}
\ss

We now describe the author's favorite model, which we describe with 
some generality, but for which very little is known, even for 
$2$-dimensional sets in $\R^3$. 

We give ourselves a finite collection of boundary pieces $\Gamma_j$, 
$0 \leq j \leq j_{max}$, (those are just closed subsets of $\R^n$ for 
the moment), and an initial competitor $E_0$ (a closed set). 
We simply want to minimize $\H^d(E)$ in the class ${\cal F}(E_0)$ 
of \underbar{sliding deformations} of $E$, which we define as follows. 
A closed set $E$ lies in ${\cal F}(E_0)$ if $E = \varphi_1(E_0)$, 
where $\{\varphi_t \}$, $0 \leq t \leq 1$, is a one-parameter family
of mappings such 
that
$$
(t,x) \to \varphi_t(x) : [0,1] \times E_0 \to \R^n
\hbox{ is continuous,}
\leqno (4.2)
$$
$$
\varphi_0(x) = x \hbox{ for } x\in E_0,
\leqno (4.3)
$$
$$
\varphi_t(x) \in \Gamma_j \hbox{ when $0 \leq j \leq j_{max}$ and } 
x\in E_0 \cap \Gamma_j,
\leqno (4.4)
$$ 
and
$$
\varphi_1 \hbox{ is Lipschitz.}
\leqno (4.5)
$$
We decided to require (4.5), mostly by tradition and to accommodate
size minimizers below, but this is negotiable. Notice that we do not 
require any quantitative Lipschitz bound for $\varphi_1$.

Of course we should check that
$$
0 < \inf_{E \in {\cal F}(E_0)} H^d(E) < +\infty,
\leqno (4.6)
$$
because otherwise the minimization problem below is not interesting.
And then we want to minimize $H^d(E)$ in the class ${\cal F}(E_0)$.

Let us just give a few comments here, starting with the
good news.

This definition seems natural for soap films. We
allow the soap film to move continuously, and we impose the
constraint that points that lie on a boundary piece $\Gamma_j$
stay on $\Gamma_j$ (but may move along $\Gamma_j$).
We also allow deformations $\varphi$ that are not injective
(including on the $\Gamma_j$); so we are allowed to pinch
and merge different parts of $E_0$. This seems to be all right 
with real soap films.

Some boundary pieces may play the role of the curve $\Gamma$
in the Plateau problems of Section 2, but we may also consider
the case when some $\Gamma_j$ are two-dimensional, and the
surface boundary is allowed to slide along $\Gamma_j$,
like a soap film that would be attached to a wall. Or we could use
a set $\Gamma_0$ that contains $E_0$ to force our films to stay
in the given region $\Gamma_0$.

The fact that we choose an initial $E_0$ also gives us some 
extra flexibility; this looks like cheating, because given a soap 
film, we can always try it as our initial set $E_0$ and see what 
happens. But at the same time, it is probable that
real soap films do something like this. And we don't have to think 
too hard about how to model each given soap film, or to wonder
about which precise topological property (for instance, belonging
to some homology class, for some group that we would need to choose) 
defines the correct class of competitors.

The definition above is not really new, even though the author did
not find more than allusions to this way of stating Plateau problems
in the literature. But Brakke's software ``Surface Evolver"
allows this as one of the main options. Maybe people just did not want 
to state a problem that looked too hard to solve.

Our definition looks like Almgren's definition of
``restricted", or $(M,\varepsilon,\delta)$-minimal sets, but in the
present situation, Almgren would tend to work on the open set
$U = \R^n \sm \big(\cup_j \Gamma_j \big)$, and would 
use the (much too strong) condition that $\varphi_t(x) = x$
for $x$ in a neighborhood of the $\Gamma_j$.
A minor difference is that Almgren typically considers only $\varphi_1$,
regardless of the existence of a one-parameter family of mappings that 
connects it to the identity; this usually makes no difference, because
in most results everything happens in a small ball contained in $U$,
and the one-parameter family can easily be obtained by convexity.
He can also play with the small parameter $\delta$, which limits
the diameter of the authorized modifications and can be used to 
forbid mappings $\varphi_1$ that are not homotopic to the identity inside $U$.

We could have taken a more restrictive approach and just required
that $\varphi_t(x) = x$ when $x\in \Gamma_j$, but the author thinks
that soap does not really act like this. For instance, for 
$2$-dimensional films in $\R^3$, with a boundary $\Gamma_j$ which
is a plane, we expect solutions that look like two half planes
with a common boundary (a line) in $\Gamma_j$, and that lie on
the same side of $\Gamma_j$. With (4.4), they should make
equal angles with $\Gamma_j$; with the stronger ``sticking boundary"
condition, any two planes should work (provided that they make an 
angle of at least $120$ degrees), and even more complicated sets
bounded by a curve in $\Gamma_j$. The definition with (4.4)
also has the advantage, since it allows more competitors, to
make it slightly easier to prove regularity results for the minimizers.

Also, we decided to define the $\varphi_t$ only on $E$, because we
do not want to force any deformation $\{ \varphi_t \}$ that the soap
may choose, to extend continuously the the whole space and respect
boundary constraints where the soap is not even present.

Of course we expect that some soap films in nature will exist just 
because they are stable local minima of $\H^d$ without really minimizing, 
even in the same class ${\cal F}(E_0)$. Even worse, the set depicted 
in Figure 6 (Section 3.e) can be deformed into a set of vanishing
measure, without even going through a set of larger Hausdorff measure, 
so the only way to exclude this deformation would be to disqualify it
because it is too long.

At the same time, we should observe that solutions of the mass-minimizing
problem of Section 2.c (where the boundary is a nice curve) provide local solutions. 
Indeed let $T$ minimize $Mass(T)$ among  solution of $\d T = S$, 
where $S$ is the current of integration on a smooth curve $\Gamma$, 
and let $E_0$ denote its support. By [HS] 
and if $\Gamma$ is $C^2$,  the set $E_0$ is regular, including along the boundary, 
and then it is easy to see that $E_0$ is locally minimal in the following sense: 
there is a small $\delta >0$ such that $\H^2(E_0) \leq \H^2(E)$ 
whenever $E$ is a deformation of $E_{0}$ in a ball $B$ of radius $\delta$.
This last means that $E \in {\cal F}(E_0)$ is obtained as in (4.2)-(4.5), 
but with mappings $\varphi_t$ such that $\varphi_t(x) = x$
for $x\in B$ and $\varphi_{t}(B) \i B$.
Of course we only get one local sliding minimizer $E_{0}$ this way, 
and it is always smooth, so the main defect of this point of view, in
the author's opinion, is that it never gives a soap film with a singularity.

Now the worse news. We are again trying to play with parameterizations
(here, with the initial set $E_0$ as a source set), and we know that it will
be hard to find optimal ones, and that we will lack compactness at
the crucial moment if we are not careful. That is, if we select
a minimizing sequence in ${\cal F}(E_0)$, the limits of convergent
subsequences will in general not lie in ${\cal F}(E_0)$.

And indeed no general existence result is known so far, even 
when $d=2$, $n=3$, there is a single piece of boundary $\Gamma$,
which is a loop, and $E_0$ is the continuous image of a disk that 
closes the loop. Notice that the Douglas solution does not help
here, because we look at a different functional which takes care
of interactions between pieces.
Thus the situation is as bad (and probably a tiny bit worse) 
than for size-minimizing currents. 

It still looks interesting to study sliding minimizers.
First because there is still a small chance that this approach 
will work is some cases, by selecting carefully a nice minimizing 
sequence before we take any limit. 
See [Da6], where a short description 
of two recent results of this type is given, but for simpler problems
where there the class ${\cal F}(E_0)$ is not given in terms boundary
pieces $\Gamma_j$ as above, but of softer topological conditions.

Also, the chances of proving existence results will probably increase 
if we understand better the regularity results for minimizers, all the
way to the boundary. We shall see in Section 7 
that such results could also be used for the solutions of the 
Reifenberg and size minimization problems (when they exist).

\ssi{\bf 4.f. Variants of the Plateau problem}
\ss

Let us just mention here that there are many other interesting 
problems where one tries to minimize $\H^d(E)$ (or some variant)
under topological conditions on $E$ (separations conditions in 
codimension $1$, homology conditions in higher codimensions, etc.)
We refer to [Da3] and [Da6] for some examples, 
but we do not develop here.

\msi{\bf 5. REGULARITY RESULTS FOR ALMGREN ALMOST MINIMAL SETS}
\ss

\ssi{\bf 5.a. Local regularity}
\ss
The following notion was introduced by Almgren, as a very good
model for studying the local regularity properties of soap films
and bubbles (among other objects). We simplify the definition
a bit, but not in a significant way.

We give ourselves an open set $U \i \R^n$, a dimension $d \leq n$,
and a nondecreasing gauge function $h : (0,+\infty) \to [0,+\infty)$, 
such that
$$
\lim_{r \to 0} h(r) = 0,
\leqno (5.1)
$$
which will measure how close we are from minimality.
Taking $h(r) = 0$ corresponds to minimality, and
$h(r) = C r^\alpha$ with $0 < \alpha \leq 1$ is a standard choice.
This is a way to accommodate slightly more complicated functionals 
than just $\H^d$.

We say that the closed subset $E$ of $U$ is a $d$-dimensional 
\underbar{almost minimal set} in $U$, with gauge function $h$, 
if the following holds.

For each closed ball $B = \overline B(x,r) \i U$, and each
Lipschitz mapping $\varphi : \R^n \to \R^n$ such that
$$
\varphi(x) = x \ \hbox{ for } x\in \R^n \sm B \ \hbox{ and }
\varphi(B) \i B,
\leqno (5.2)
$$
then
$$
\H^d(E\cap B) \leq \H^d(\varphi(E)\cap B) + r^d h(r);
\leqno (5.3)
$$
we also demand that
$$
\H^d(E\cap B) < +\infty
\leqno (5.4)
$$
for all $B \i U$ as above, to avoid stupidly large sets.

Notice again that $\varphi$ is allowed not to be injective,
and that since $B$ is convex, we could easily connect
$\varphi$ with the identity by a one parameter family of 
mappings that satisfy (5.2)
(take $\varphi_t(x) = t \varphi(x) + (1-t) x$).

It is  clear that any solution of the sliding Plateau problem 
mentioned in Section 4.d is a minimal set in the complement of
the $\d \Gamma_j$. (We prefer to say $\d \Gamma_j$, because some
$\Gamma_j$ may be $n$-dimensional and designed to contain
$E_0$ and all its competitors.) 
We shall see in Section 7 that this also applies to minimizers of 
Reifenberg's problem, or supports of size minimizers.

The local regularity in $U$ of the almost minimal sets 
was started by Almgren [Al4], 
and continued in [Ta2], [DS], 
[Da1,4,5], and others.  
For general dimensions and codimensions, we get that,
modulo a set of vanishing $\H^d$-measure, the almost minimal set
$E$ is locally Ahlfors-regular, rectifiable, and even uniformly 
rectifiable with big pieces of Lipschitz graphs. When
$E$ is $2$-dimensional in $\R^3$, Taylor [Ta2] 
proved that it is locally $C^1$-equivalent to a minimal
cone (a plane, or a set $Y$ or $T$ as in Section 3.d), 
and this was partially extended to higher ambient dimensions
in [Da4,5]. 
Also see [Da6] for a slightly 
more detailed survey of these results.

The reader may regret that this is not very smooth, but
notice that the $C^1$ description of J. Taylor is nearly optimal,
since after all the minimal cones above are almost minimal sets;
the situation for $d > 2$ is widely open.

A last property of the notion of almost minimal sets that may turn out
to be useful is its stability under limits [Da1]: 
if $\{ E_k \}$ is a sequence of almost minimal sets, 
with the same gauge function $h$, and that converges 
(relative to the Hausdorff distance, and after
cutting out unneeded sets of vanishing measure) to a limit $E$,
then $E$ is also almost minimal, with the same gauge function $h$.
In addition, 
$$
\H^d(E\cap V) \leq \liminf_{k \to +\infty} \H^d(E_k\cap V)
\leqno (5.5)
$$
for every relatively compact open subset $V$ of $U$; that is,
the restriction of $\H^d_{\vert V}$ to our sequence is 
lowersemicontinuous.

\ssi{\bf 5.b. Regularity near the boundary}
\ss

Not so many regularity results exist that go all the way to 
the boundary. The author knows about [All] 
(for varifold, near a flat point), 
and [HS] 
and [Wh] (for mass minimizers). 
Also see [LM2], Figure 5.3, 
or [Mo5], Figure 13.9.3 on page 137 in my third edition, 
for a conjecture about the 
types of singularities of a soap film near the boundary, and
[Br3] for a description of soap films 
near a point where they leave a boundary curve.

It is nonetheless interesting to see what can be done near the 
boundary, in the context of Almgren almost minimal sets.
Let the closed boundaries pieces $\Gamma_j$, $0 \leq j \leq j_{max}$, 
be chosen, as in Section 4.e, and let $E \i \R^n$ be a closed set, 
with $\H^d(E) < + \infty$.

\ms\proclaim Definition 5.6.
We say that $E$ is \underbar{almost minimal, with
the sliding conditions}
defined by the boundaries $\Gamma_j$, if 
$$
\H^d(E\cap B) \leq \H^d(\varphi_1(E)\cap B) + r^d h(r)
\leqno (5.7)
$$
(as in (5.3), but) whenever $B = \overline B(x,r) \i \R^n$ and the
$\{\varphi_t \} : E \to \R^n$, $0 \leq t \leq 1$, are such that
$$
(t,x) \to \varphi_t(x) : [0,1] \times E \to \R^n
\hbox{ is continuous,}
\leqno (5.8)
$$
$$
\varphi_0(x) = x \hbox{ for } x\in E,
\leqno (5.9)
$$
$$
\varphi_t(x) \in \Gamma_j \hbox{ when $0 \leq j \leq j_{max}$ and } 
x\in E \cap \Gamma_j,
\leqno (5.10)
$$ 
$$
\varphi_1 \hbox{ is Lipschitz.}
\leqno (5.11)
$$
as in (4.2)-(4.5), and, for $0 \leq t \leq 1$,
$$
\varphi_t(x) = x \ \hbox{ for } x\in E \sm B \hbox{ and } 
\varphi_t(B) \i B.
\leqno (5.12)
$$

\ms
Thus, when we take $h = 0$, we get that $E$ is a minimizer
if it solves the problem of Section 4.e with $E_0 = E$.
The local almost minimality property above amounts to a little
less than taking infinitely many boundary pieces, equal to the 
points of $\R^n \sm U$, and applying the definition with (5.7)-(5.12).

Notice that the fact that we are now allowed to move points of the
boundaries is good for us, because it allows more competitors and
we can hope to get some direct information at the boundary.

Anyway, the author decided to run all the local regularity
proofs he knows, and try to extend them to the boundary. 
So far (subject to additional proofreading), local Ahlfors-regularity, 
rectifiability, and limiting theorems seem to be under control, 
assuming for instance that the boundary pieces are all faces 
(possibly of all dimensions) of a dyadic grid, or of the image of a 
dyadic grid by a $C^1$ diffeomorphism. 
The next stages for an extension of J. Taylor's result would be to get 
some control on the blow-up limits of $E$ at a boundary point
(a new, larger list of cones that we also need to determine), 
and some version of the monotonicity 
of density, including for balls that are centered near the boundary,
but not exactly on the boundary. This last could be problematic (but
worth trying).

Quite a few interesting things can happen at the boundary,
even when $d=2$ and the boundary is a smooth curve.
See Figure 7 and other pictures of soap films 
over the Borromean Rings in K. Brakke's website, 
and a description [Br3] 
of what may happen when the minimal surface leaves the curve 
(as in the examples of Figures 2.d and 7). 

\msi{\bf 6. AMNESIC SOLUTIONS OF THE PLATEAU PROBLEM}
\ss
In this short section we try to clarify some issue about
existence results for Plateau problems. To make things simpler,
let us consider the case of sliding minimizers, with only one piece 
of boundary, a smooth closed curve $\Gamma \i \R^3$. We 
parameterize $\Gamma$ by $\gamma : \d D \to \Gamma$,
extend $\gamma$ to the unit disk, and get an initial set $E_0$.
Then we consider the set ${\cal F}(E_0)$ of deformations 
$E = \varphi_1(E_0)$ of $E_0$, where the $\varphi_t$ preserve $\Gamma$
(see near (4.2)). We would like to find $E \in {\cal F}(E_0)$ such that
$\H^d(E) = m$, where
$$
m = \inf\big\{ \H^2(F) \, ; F \in {\cal F}(E_0) \big\}.
\leqno (6.1)
$$
Notice that we look for an absolute minimizer here, which in principle should
make things simpler, but at the same time would forbid us from using 
solutions of other problems in other classes, such as mass minimizers.

An approach like the one we shall describe here would also make sense in 
more general situations, or for the Reifenberg problems of Section 2.b, but
we shall not elaborate. Also, the reader may want to skip part of
the construction of an amnesic solution just below, and
go directly to one of the main points of the section, which is
what the author means by amnesic solutions and why he thinks they are
not entirely satisfactory.

Select a minimizing sequence $\{ E_k \}$. That is,
$E_k \in {\cal F}(E_0)$ and $\H^2(E_k) \leq m + \varepsilon_k$,
where $\varepsilon_k$ tends to $0$ (and we may assume that
$\varepsilon_k \leq 2^{-k}$). Let $B$ be a large closed ball 
that contains $\Gamma$ and $E_0$; we can assume that every $E_k$ 
is contained in $B$ (otherwise, project radially on $B$ and 
you will get a competitor which is at least as good). 
Then replace $\{ E_k \}$ with some subsequence that tends 
to a limit $E$ (for the Hausdorff distance on compact sets). 

We would like $E$ to be a minimizer for our problem, but if
we don't pay attention, this will surely not happen. Indeed,
we can easily choose $E_k$ with lots of long and thin hair,
with almost no area, but so that $E_k$ is $2^{-k}$-dense in $B$.
If we do this, we get $E=B$, which is not even be $2$-dimensional.

But there is a way to pick the $E_k$ carefully, so that the limit
$E$ represents a more respectable attempt. In [Re1], 
Reifenberg does this by carefully cutting the hair of an initial
$E_k$. Let us say two words about a slightly different way; we refer
to [Da6] 
for a more detailed account of the strategy and some applications.

We start from our initial $E_k\in {\cal F}(E_0)$
and first use a construction of Feuvrier [Fv2] 
to build a sort of dyadic grid ${\cal G}_k$ adapted to $E_k$. 
That is, instead of decomposing $\R^3$ into the almost disjoint union
of dyadic cubes of size $2^{-l}$ for some very small $l$, we use
polyhedra instead of cubes. The polyhedra will all have a diameter
comparable to some $2^{-l}$ (the mesh of ${\cal G}_k$), and 
we make sure that $2^{-l} << 2^{-k}$ in the construction.
We do not have a lower bound for the mesh, but we shall not need 
one.

We want to replace $E_k$ with a Federer-Fleming projection on the
grid ${\cal G}_k$, but since we are afraid of $\Gamma$, we shall
only do this reasonably far from $\Gamma$. Denote by $V_k$
the union of all the polyhedra of ${\cal G}_k$ that lie at
distance at least $2^{-k}$ from $\Gamma$. On $V_k$, we 
replace $E_k$ with its Federer-Fleming projection on the grid.
That is, for each polyhedron $Q \i V_k$, we take an interior point
$x_Q \in Q \sm E_k$, and replace $E_k \cap Q$ with its projection 
on the boundary $\d Q$, where the projection is the radial projection 
centered at $x_Q$. Outside of $V_k$, we change nothing. This gives a 
new set $E'_k$. 

Then we do the same construction in the $2$-dimensional 
faces of the grid. Just consider the faces $F$ that are contained 
in the interior of $V_k$, and such $E'_k$ does not fill the 
interior of $F$. That is, such that we can find an interior point 
$x_F \in F \sm E'_k$; we then use this point to
project $E'_k \cap F$ to the boundary of $F$. These manipulations are
independent, and we don't need to define the projections elsewhere,
because now $E'_k \cap V_k$ is contained in the union of the $2$-faces.
We do this projection for all the faces $F$ where this is
possible, and get a new set $E^{(2)}_k$. Notice that
$E^{(2)}_k\in {\cal F}(E_0)$, because $E_k\in {\cal F}(E_0)$,
the Federer-Fleming projections are deformations, and we made sure
not to move anything near $\Gamma$.

By the construction of the adapted grid (and maybe by choosing
a little more carefully, by a Fubini argument, the place $\d V_k$ 
where we do the interface between the identity and the Federer-Fleming
projections), we get that 
$E^{(2)}_k \leq \H^d(E_k) +2^{-k} \leq m + 2^{-k+1}$; 
the general point is that we construct ${\cal G}_k$
so that $E_k$ is often very close to $2$-faces of the grid,
so that the projections will not make $\H^2(E_k)$ much larger.

At last $E^{(2)}_k$ looks nicer inside $V_k$. That is,
let $V'_k \i V_k$ denote the union of the polyhedra $Q$ of ${\cal G}_k$
such that every polyhedron $R$ of the grid that touches $Q$ is 
contained in $V_k$ (we just remove something like one exterior layer of 
polyhedra). Then inside $V'_k$, $E^{(2)}_k$ is composed of entire 
$2$-faces of the grid ${\cal G}_k$, plus possibly parts of $1$-faces 
of ${\cal G}_k$. 
We could continue one more step to get rid of the $1$-faces that are not 
entirely covered, but this will not be needed.

Now $E^{(2)}_k$ is not yet our cleaner competitor. 
We look for deformations $F$ of $E^{(2)}_k$ inside $V'_k$, that are 
also composed, inside $V'_k$, of entire $2$-faces of the grid, 
plus possibly parts of $1$-faces of the grid, and for which 
$\H^2(F)$ is minimal. Here deformation of $E^{(2)}_k$ inside $V'_k$
means that $F = \varphi_1(E^{(2)}_k)$, where the $\varphi_t$,
$0 \leq t \leq 1$, are such that $\varphi_t(x) = x$ for $x\in \R^3 
\sm V'_k$ and $\varphi_t(V'_k) \i V'_k$, in addition to the usual
(5.8), (5.9), and (5.11). Minimizers for this problem exist
trivially, because modulo sets of dimension 1, there is only 
a finite number of sets $F$.
We select a minimizer $F_k$, which is our cleaner competitor.

A second property of the adapted grids is that we have uniform lower 
bounds on the angles of the faces of the polyhedra that compose them.
Because of this, the set $F_k$ has some regularity inside $V'_k$:
its minimality property implies that it is ``Almgren quasiminimal"
far from $\R^3 \sm V'_k$, and this is enough to imply some 
lowersemicontinuity of $\H^d$, restricted to the sequence $\{ F_k \}$.
That is, (5.5) holds for the $F_k$, and for any open set $V$
which is compactly contained in $\R^3 \sm \Gamma$.

Here we are going a little fast; all these things need proofs,
but we just want to give an idea of why it helps to make this complicated 
construction. We should also say that our restriction to $n=3$ and
$d=2$ is not needed for this part of the argument; it just makes
things more explicit.
In other contexts (currents, varifolds), this stage could be less
painful, because there exist powerful lowersemicontinuity results
that can take care of the analogue of (5.5).

We now take a subsequence so that $\{ F_k \}$ converges to a set $F$,
and we would like to say that $F$ is a solution of our problem.
First, because of (5.5) (and because $\H^2(\Gamma) = 0$), we really get
that 
$$
\H^2(F) \leq \liminf_{k \to +\infty}  \H^2(F_k) 
\leq \liminf_{k \to +\infty}  \H^2(E^{(2)}_k) \leq m,
\leqno (6.2)
$$
so $F$ looks like a good candidate. We can also say more about
$F$ in $U = \R^3 \sm \Gamma$. From the fact that $F_k$ nearly
minimizes $\H^2$ in the class ${\cal F}(E_0)$, the fact that
${\cal F}(E_0)$ is stable under local deformations in compact balls
$B \i U$ (as in (5.2)),
and the lowersemicontinuity property (5.5), one can deduce that
$$
F \hbox{ is an Almgren minimal set in $U = \R^3 \sm \Gamma$.}
\leqno (6.3)
$$
That is, (5.3) holds, with $h(r) = 0$.

\ms
This is not so bad, but let us now say why we may call $F$
an amnesic solution. In order to get a real solution to the initial
problem, we would still need to check that $F \in {\cal F}(E_0)$,
and this would involve finding a deformation $\{ \varphi_t \}$
from $E_0$ to $F$. This will be hard if we keep the argument
as it is suggested, because we have no control on the various
deformations $\{ \varphi_t^k \}$ associated to the $F_k$. 
The situation is not completely hopeless, and one can prove some
existence results in similar contexts (but not yet for Plateau's 
problem with a curve);
see [Da6].  

Return to our set $F$. It looks like a solution of Plateau's
problem, especially far from $\Gamma$; in addition to (6.3), something
probably remains from the initial problem that was used to construct 
$F$. We know that $F$ lies in the closure of ${\cal F}(E_0)$, 
and in some problems this may be enough information, but otherwise
it is hard to say exactly what properties are preserved. 
In the specific sliding problem here, an intermediate information would be 
to prove that $F$ is a sliding minimal set, i.e., a solution of the 
problem above, where we replace $E_0$ with $F$ itself, but even then
we will not be entirely happy, because some information may have been 
lost in the limit.

Part of the reason for this section was to insist on the difference,
in the author's view, between a complete solution of a Plateau problem
(like the one above, or some of the more classical variants of 
Section 2) and an amnesic solution, which is often easier to produce.
The difference may be subtle though, because it may turn out that in 
fact all the amnesic solutions are true solutions. That is, it seems 
hard to find an example of an amnesic solution with some clearly 
missing feature.

\msi{\bf 7. REIFENBERG MINIMIZERS AND SUPPORTS OF SIZE MINIMIZERS
ARE SLIDING MINIMAL}
\ss
In this last section we record the fact that in many cases, the solutions of
generalizations of the Reifenberg problem of Section 2.b, as well as
the size minimizing current problem of Section 2.c,
give sliding minimal sets. Thus boundary regularity results for 
sliding minimal sets may be used in these contexts as well, even 
though one may object that they may be easier to obtain directly.

Also notice that just by restricting to the complement $U$ 
of the boundary set, we see that locally in $U$, these problems 
yield locally Almgren minimal sets, for which the interior 
regularity results of Section 5.a apply.
This fact is apparently part of the folklore.

We start with variants of the Reifenberg problem.
Choose integers $0 < d < n$, and a notion of homology 
(including the choice of an abelian group $G$). Any usual
choice will do; we shall just use the homotopy invariance
and the fact that $\d$ is a natural map.

Also choose a closed boundary set $\Gamma \i \R^n$,
and a collection $\{ \gamma_j \}$, $j\in J$, of elements 
of the homology group $H_{d-1}(\Gamma)$. Let us 
assume that
$$
\H^d(\Gamma) = 0;
\leqno (7.1)
$$
this will allow us to add $\Gamma$ to our competitors
$F$ at no cost, and the definitions will be much easier to apply.

Then let $\cal F$ denote the class of closed sets 
$F \i \R^3$ that contain $\Gamma$ and for which 
$i_\ast(\gamma_j) = 0$ in 
$H_{d-1}(F)$ for all $j\in J$, where
$i_\ast : H_{d-1}(\Gamma) \to H_{d-1}(F)$ 
is the natural map coming from the injection 
$i : \Gamma \to F$. 

We shall say that $E$ is Reifenberg-minimal
(relative to the choices above) if $E \in {\cal F}$
and
$$
\H^d(E) = \inf_{F \in {\cal F}} \H^d(F) < +\infty.
\leqno (7.2)
$$

\proclaim Proposition 7.3.
If (7.1) holds and $E$ is Reifenberg-minimal, then it is also
a minimal set of dimension $d$, as in Definition 5.6, with the 
sliding conditions defined by the unique boundary piece $\Gamma$. 

\ms
Indeed, let $B$ and the $\varphi_t$ be as in Definition 5.6; 
we want to show, as in (5.7), that 
$\H^d(E\cap B) \leq \H^d(\varphi_1(E)\cap B)$. 
Since here $\H^d(E) < +\infty$, we can add $\H^d(E\sm B)$
to both terms, and this simplifies
to $\H^d(E) \leq \H^d(\varphi_1(E))$.
Then set $F = \Gamma \cup \varphi_1(E)$; by (7.1), it is
still enough to prove that $\H^d(E) \leq \H^d(F)$,
and for this (and by (7.2)) we just need to prove that $F$ 
lies in $\cal F$.

So let $j\in J$ be given, and let us check that
$i_{\ast}(\gamma_j) = 0$ in $H_{d-1}(F)$, where 
$i_{\ast} : H_{d-1}(\Gamma) \to H_{d-1}(F)$ 
comes from the injection $i : \Gamma \to F$. 
Denote by $i_0 : \Gamma \to E$
the initial injection. Since $E \in {\cal F}$, 
$i_{0,\ast}(\gamma_j) = 0$ in $H_{d-1}(E)$, 
which means that there exists a chain $\sigma$ in $E$
such that $\d \sigma = i_{0,\ast}(\gamma_j)$.

Apply the mapping $\varphi_1$ to this. 
Set $\sigma_1 = \varphi_{1,\ast}(\sigma)$; this
is a chain in $\varphi_1(E)$, and
$$
\d \sigma_1 = \d \varphi_{1,\ast}(\sigma)
= \varphi_{1,\ast}(\d \sigma)
= \varphi_{1,\ast}(i_{0,\ast}(\gamma_j))
= (\varphi_1 \circ i_0)_\ast(\gamma_j)
\leqno (7.4)
$$
because $\d$ is natural and by definition of $\sigma$.

Notice that for $0 \leq t \leq 1$, $\varphi_t$ is defined 
on $\Gamma$ (because $E$ contains $\Gamma$ since $E\in \cal F$),
and $\varphi_t(\Gamma) = \varphi_t(\Gamma \cap E) \i \Gamma$
by (5.10). Call $\varphi : \Gamma \to \Gamma$ the restriction of 
$\varphi_1$ to $\Gamma$ (to distinguish), and set
$$
\gamma'_{j} = \varphi_\ast(\gamma_j)\in H_{d-1}(\Gamma).
\leqno (7.5)
$$
Recall that $\varphi_0(x) = x$ on $\Gamma$ (by (5.9)),
so the $\varphi_t$, $0 \leq t \leq 1$, provide a homotopy
from the identity to $\varphi$ (among continuous mappings from $\Gamma$
to $\Gamma$). Now we use the invariance of homology under homotopies,
and get that $\gamma'_{j} = \gamma_j$ in $H_{d-1}(\Gamma)$ because
$\varphi$ and the identity induce the same mapping on $H_{d-1}(\Gamma)$.

Next denote by $i_1 : \Gamma \to \varphi_1(E)$ the inclusion,
and notice that $\varphi_1 \circ i_0 = i_1 \circ \varphi$ 
(as a map from $\Gamma$ to $\varphi_1(E)$). Then (7.4) yields
$$
\d \sigma_1 = (\varphi_1 \circ i_0)_\ast(\gamma_j)
= (i_1 \circ \varphi)_\ast(\gamma_j)
= i_{1,\ast}(\gamma'_{j})
= i_{1,\ast}(\gamma_{j})
\leqno (7.6)
$$
in $H_{d-1}(\varphi_1(E))$; hence $i_{1,\ast}(\gamma_{j}) = 0$
in $H_{d-1}(\varphi_1(E))$, because it is a boundary. 

We compose with a last injection $i' : \varphi_1(E) \to F$, 
and get that
$i_\ast(\gamma_j) = (i' \circ i_1)_\ast(\gamma_j)
= i'_\ast(i_{1,\ast}(\gamma_{j})) = 0$
in $H_{d-1}(F)$, as needed.

So $F \in {\cal F}$, and $E$ is a minimal set, as promised.
The main goal of the detailed computations above was to
convince homology nincompoops such as the author that nothing 
goes wrong with the arrows. 
\qed

\msi{\bf Remark 7.7.}
Our assumption (7.1) is not infinitely shocking
(for instance, curves in $\R^n$ satisfy this when $d=2$,
and (7.1) still allows lots of room for chains of dimension
$d-1$), but it is not clear that it should really be there.

In the definition of $F \in \cal F$, it makes sense
to take $F$ closed, not necessarily containing $\Gamma$,
and to say that $i_\ast(\gamma_j) = 0$ in $H_{d-1}(E \cup \Gamma)$,
where $i$ now denotes the injection from $\Gamma$ to $E \cup \Gamma$.

But if $\H^d(\Gamma) > 0$, we have some problems with
the proof above. We used the fact that $E$ contains $\Gamma$
to say that our $\varphi_t$ are defined on $\Gamma$
(and not just $E \cap \Gamma$), and this was used to compute
$i_\ast(\gamma_j)$. We don't want to include $\Gamma$ if
$H^d(\Gamma) > 0$, because $E \cap \Gamma$ should not be minimal
(maybe a big part of $\Gamma$ is just useless).

If for some reason we know that our mapping $(t,x) \to \varphi_t(x)$,
from $[0,1] \times (E\cap \Gamma)$ to $\Gamma$,
has a continuous extension from $[0,1] \times\Gamma$ to $\Gamma$,
we can compute as above.

Otherwise, and if we can choose chains in $\Gamma$ that represent 
the $\gamma_j$, and whose supports are all contained in a close set 
$\Gamma' \i \Gamma$, with $\H^d(\Gamma') = 0$, we can also try to 
finesse the issue by showing that $E\cup \Gamma'$ is minimal, 
with the sliding condition associated to $\Gamma$, 
but this is hard if the support of $\sigma$
is not contained in $E \cup \Sigma'$ (but has a significant 
piece in the rest of $\Gamma$). Let us not try to get a statement.

\ms
Next we turn to size-minimizing currents, as in Section 2.c,
but again we shall consider a slightly more general problem.

Let $0 < d < n$ be integers, and let $\Gamma$ be a compact subset
of $\R^n$. Also let $S$ be an integral current of dimension $d-1$ 
in $\R^n$, with $\d S = 0$ and with support in $\Gamma$. 
This last just means that $\langle S,\omega \rangle = 0$ 
when the $(d-1)$-form $\omega$ is supported in $\R^n \sm \Gamma$. 

Next let $\cal S$ denote the collection of all the integral current $S'$
of dimension $d-1$ that are supported on $\Gamma$ and homologous
to $S$ on $\Gamma$. This last means that there exists a $d$-dimensional 
current $V$, supported on $\Gamma$, 
and such that $\d V = S-S'$.

Finally denote by $\cal T$ the set of integral currents $T$
such that $\d T \in {\cal S}$. This is our set of competitors,
and we want to minimize $Size(T)$ over $\cal T$.

Notice that this setting includes the simple example of Section 2.
Indeed, let $\Gamma$ be a smooth orientable surface of dimension 
$d-1$, and let $S$ be the current of integration on $\Sigma$.
Notice that ${\cal S} = \{ S \}$, because if $V$ is a $d$-dimensional 
current supported on $\Gamma$, then $V = 0$. In this case we just
want to minimize $Size(T)$ over the solutions of $\d T = S$, as above.

But we could also take for $\Gamma$ some $2$-dimensional torus 
along a closed curve, pick a closed loop on $\Gamma$, let $S$ be the
current of integration over that loop, and $\cal S$ will allow the 
current of integration over any Lipschitz deformation of that loop
in $\Gamma$ (see the proof below). This is a way to encode a possibly
sliding boundary.

\ms\proclaim Proposition 7.8.
Let $\Gamma$, $\cal S$, and $\cal T$ be as above, 
and let $T\in \cal T$ be such that
$$
Size(T) = \inf_{R \in {\cal T}} Size(R) < +\infty.
\leqno (7.9)
$$
Denote by $Z$ the closed support of $\d T$, and assume that
$$
\H^d(Z) = 0.
\leqno (7.10)
$$
Also assume that $\Gamma$ is a Lipschitz neighborhood retract.
Let $E$ denote the support of $T$; then $E \cup Z$ 
is a minimal set of dimension $d$, as in Definition 5.6, with the 
sliding conditions defined by the single boundary piece $\Gamma$.

\ms
Probably we can prove (7.10) (instead of assuming it) in some cases, 
but let us not bother.
Our conclusion that $E \cup Z$ is minimal, rather than just $E$, 
looks unpleasant (the Hausdorff measures are the same, but there is 
a difference because the boundary constraint (5.10)
for the competitors may be different). Under fairly weak assumptions,
one can show that if $F$ is minimal (with some sliding boundary 
conditions), then the closed support of the restriction of $\H^d$
to $F$ is minimal too (with the same sliding boundary 
conditions), which could be simpler. 

Our neighborhood retract assumption is just that there is a 
small neighborhood $W$ of $\Gamma$, and a mapping
$h : W \to \Gamma$, which is Lipschitz and such that $h(x)=x$
for $x\in \Gamma$.
The reader should not pay too much attention to all these details,
which are mostly technical.

So let $T$ be a size minimizer, as in the statement. 
Recall that $T$ has an expression like (2.6), namely, 
$$
\langle T,\omega \rangle = \int_A m(x) \; \omega(x)\cdot\tau(x) 
\, d\H^d(x),
\leqno (7.11)
$$
when $\omega$ is a (smooth) $d$-form,  where $A$ is a rectifiable set,
$\tau(x)$ is a $d$-vector that spans the approximate tangent $d$-plane 
to $A$ at $x$, and $m$ is an integer-valued multiplicity function on $A$, 
integrable against ${\bf 1}_E d\H^d$. Also recall that 
$$
Size(T) = \H^d(A'), \hbox{ with }
A' = \big\{ x\in A \, ; \, m(x) \neq 0 
\big\}.
\leqno (7.12)
$$
The Borel support $A'$ may be strictly smaller than $E$, which 
is in fact equal to the closure of $A'$, so we need to be careful 
about $\H^d(E\sm A')$. We claim that
$$
\H^d(E\sm (Z \cup A')) = 0.
\leqno (7.13)
$$
That is, away from the support $Z$  of $\d T$,
$E \sm A'$ has vanishing Hausdorff measure. Let merely
sketch the proof here. 
There is a monotonicity formula, 
which says that for $x\in E \sm Z$, the density function
$r \to r^{-d} \H^d(A' \cap B(x,r))$ is nondecreasing
on $(0,\dist(x,Z))$. This is, classically, obtained by replacing $T$
on any small ball $B(x,r)$ with the cone over its slice on $\d B(x,r)$,
comparing, and then integrating the result. 
Next let $\omega_d$ denote the $\H^d$-measure
of the unit ball in $\R^d$; since $A'$ is rectifiable,
$$
\lim_{r \to 0} r^{-d} \H^d(A' \cap B(x,r)) = \omega_d
\ \hbox{ for $\H^d$-almost every } x\in A';
\leqno (7.14)
$$
see [Ma], for instance, 
for this and the next density results. 
Hence (by monotonicity)
$$
\H^d(A' \cap B(x,r)) \geq \omega_d r^d
\ \hbox{ for } 0 < r < \dist(x,Z)
\leqno (7.15)
$$
for almost-every $x\in A'$, and hence (take a limit) for every $x\in E$ too.
But by a standard density theorem, 
$$
\lim_{r \to 0} \H^d(A' \cap B(x,r)) = 0
\ \hbox{ for $\H^d$-almost every }
x\in \R^n \sm A', 
\leqno (7.16)
$$
and (7.13) follows by comparing (7.16) and (7.15). 

\ms
We are now ready to check that $E' = E \cup Z$ is minimal,
with sliding conditions associated to $\Gamma$. 
Let $B$ and the $\varphi_t$ be as in Definition 5.6 (for $E'$); 
we want to show that $\H^d(E' \cap B) \leq \H^d(\varphi_1(E') \cap B)$,
and since 
$$\eqalign{
\H^d(E') &= \H^d(E\cup Z) = \H^d(E\sm Z) + \H^d(Z) = \H^d(E\sm Z) 
\cr&
\leq \H^d(E\sm (Z \cup A')) + \H^d(A')
= \H^d(A') < +\infty
}\leqno (7.17)
$$
by (7.10), (7.13), and (7.12), we can add or subtract 
$\H^d(E'\sm B) = \H^d(\varphi_1(E')\sm B)$, and so
we just need to check that 
$$
\H^d(E') \leq \H^d(\varphi_1(E')).
\leqno (7.18)
$$
We want to build a competitor for $T$, and logically we 
use $\varphi = \varphi_1$ to push $T$ forward 
and set $T_1 = \varphi_\sharp T$. Recall (for instance from 
the end of 4.1.7 in [Fe1]) 
that when $\varphi$ is smooth, $\varphi_\sharp T$ is defined by
$\langle \varphi_\sharp T,\omega \rangle 
= \langle T,\varphi^\sharp \omega \rangle$, for every $d$-form 
$\omega$, where at each point $x$, $(\varphi^\sharp\omega)(x)$ is obtained 
by applying to $\omega(x)$ the ($d$-linear version of the) differential of 
$\varphi$ at $x$. Since $T$ is given by (7.11), we get that
$$
\langle T_1,\omega \rangle = \langle T,\varphi^\sharp \omega \rangle
= \int_A m(x) \; (\varphi^\sharp\omega)(x)\cdot\tau(x) 
\, d\H^d(x),
\leqno (7.19)
$$
which we may transform into an integral on $\varphi(A)$,
with a multiplicity $m(y)$ that is a sum of multiplicities
$m(x)$, $x\in \varphi^{-1}(y)$.
For general Lipschitz functions $\varphi$, 
$ \varphi_\sharp T$ is defined by a limiting argument
(see 4.1.14 in [Fe1]), 
but it is not hard to see that 
$T_1 = \varphi_\sharp T$ is an integrable current,
associated to the rectifiable set $A_1 = \varphi_1(A)$,
and whose Borel support $A'_1$ is contained 
in $\varphi_1(A')$; it may be strictly smaller because
two different pieces of $A'$ may be sent to the same
piece of $\varphi_1(A')$, and the multiplicities may
end up canceling. So
$$
Size(T_1) \leq \H^d(\varphi_1(A')).
\leqno (7.20)
$$
Our next task is to show that $T_1 \in {\cal T}$, because
as soon as we do this, the minimality of $T$ will yield
$$
\H^d(E') \leq \H^d(A') = Size(T) 
\leq Size(T_1) \leq \H^d(\varphi_1(A'))
\leq \H^d(\varphi_1(E'))
\leqno (7.21)
$$
by (7.17) and because $A' \i E'$; then (7.18) and the conclusion will 
follow.

We know from 4.1.14 in [Fe1] 
that $\d T_1 = \d(\varphi_\sharp T)  =\varphi_\sharp (\d T)$.
Since $T_1$ is an integral current, we just need to show that
$\d T_1 \in { \cal S}$. The support of $\d T_1$ is contained 
in $\varphi(Z)$. Let us check that
$$
Z \i  (E \cup Z) \cap \Gamma = E' \cap \Gamma 
\ \hbox{ and } \ 
\varphi(Z) \i \Gamma;
\leqno (7.22)
$$
indeed $Z$, the closed support of $\d T$, is contained in $\Gamma$ 
because $\d T \in {\cal S}$; the first part now follows from 
the definition of $E'$, and then
$\varphi(Z) = \varphi_1(Z) \i \Gamma$ by (5.10) for $E'$.

So $\d T_1$ is supported in $\Gamma$. 
We still need to show that $\d T_1$ is homologous to $S$
on $\Gamma$. Since $T \in \cal T$, $\d T$ is homologous
to $S$, and it is enough to show that $\d T_1$ is homologous
to $\d T$.

We would like to use the homotopy, from the origin to $\varphi = \varphi_1$, 
given by the $\varphi_t$, but there is a minor difficulty,
because we did not assume $\varphi_t(x)$ to be a Lipschitz function 
of $x$ and $t$. This is where we shall need a smoothing argument
and our assumption that $\Gamma$ is a Lipschitz neighborhood retract.

Recall that $\varphi_t$ is defined on $Z$, 
because $Z \i E' \cap \Gamma$ (by (7.22)), and in addition
$\varphi_t(Z) \i \Gamma$, again by (5.10) as in (7.22).

Let $\varepsilon > 0$ be so small that
$W$ contains a $2\varepsilon$-neighborhood of $\Gamma$, and
let $\psi : [0,1] \times Z \to \R^n$ be a smooth function
such that $|\psi(t,x)-\varphi_t(x)| \leq \varepsilon$ for 
$(t,x) \in [0,1] \times Z$.
Such a function is easy to obtain: first extend $(x,t) \to \varphi_t$,
and then smooth it.
Next define a homotopy $\{ \psi_t \}$, $0 \leq t \leq 3$,
from the identity (on $Z$) to $\varphi$ by
$$\eqalign{
\psi_t(x) &= (1-t) x + t \psi(0,x) \hskip 1.5cm \hbox{ for } 0 \leq t \leq 1,
\cr
\psi_t(x) &= \psi(t-1,x) \hskip 2.7cm \hbox{ for } 1 \leq t \leq 2,
\cr
\psi_t(x) &= (t-2) \varphi(x) + (3-t)\psi(1,x) \hskip 0.1cm \hbox{ for } 2 \leq 
t \leq 3.
}\leqno (7.23)
$$
Recall that $\varphi_t(x) \in \Gamma$ for $x\in Z$, 
so $\dist(\psi_t(x),\Gamma) \leq \varepsilon$ for 
$(t,x)\in [0,3] \times Z$, and we can define $\xi(t,x) = \xi_t(x) = h(\psi_t(x))$
for $(t,x)\in [0,3] \times Z$, where $h : W \to \Gamma$ is the Lipschitz 
retraction of the statement, and this defines a Lipschitz homotopy from the 
identity to $\varphi$, with values in $\Gamma$. Now the homotopy 
formula for currents (4.1.9 and 4.1.14 in [Fe1]) 
says that
$$
T_1 - T = \varphi_\sharp T - T = \d \xi_\sharp([0,3] \times T) 
+ \xi_\sharp([0,3] \times \d T);
\leqno (7.24)
$$
we take the boundary and get that
$$
\d T_1 - \d T = \d \xi_\sharp([0,3] \times \d T),
\leqno (7.25)
$$
which is fine because $\xi_\sharp([0,3] \times \d T)$
is a current of degree $d$ supported in $\Gamma$. Thus
$\d T_1$ is homologous to $\d T$ on $\Gamma$,
$T \in \cal T$, and this completes our proof of Proposition 7.8.
\qed

\bigskip                             
\centerline {REFERENCES }
\ss

\item {[All]} W. K. Allard,
On the first variation of a varifold. 
Ann. of Math. 95 (1972), 417--491.
\item {[Al1]} F. J. Almgren,
Plateau's problem: An invitation to varifold geometry.
W. A. Benjamin, Inc., New York-Amsterdam 1966, xii+74 pp.
\item {[Al2]} F. J. Almgren, Some interior regularity theorems 
for minimal surfaces and an extension of Bernstein's Theorem.
Ann. of Math. (2), Vol. 84 (1966), 277--292. 
\item {[Al3]} F. J. Almgren, Existence and regularity 
almost everywhere of solutions to elliptic variational problems 
among surfaces of varying topological type and singularity structure.
Ann. of Math. (2) 87 (1968) 321--391. 
\item {[Al4]} F. J. Almgren, Existence and regularity almost everywhere 
of solutions to elliptic variational problems with constraints.
Memoirs of the Amer. Math. Soc. 165, volume 4 (1976), i--199.
\item {[Al5]} F. J. Almgren, 
$Q$-valued functions minimizing Dirichlet's integral 
and the regularity of area minimizing rectifiable currents up to 
codimension two.  
Bull. Amer. Math. Soc. (N.S.) 8 (1983), no. 2, 327Ð328.
\item {[Br1]} K. Brakke, Minimal cones on hypercubes.    
J. Geom. Anal. vol. 1 (1991), 329--338. 
\item {[Br2]} K. Brakke, The surface evolver.
Experiment. Math. 1 (1992), 141--165. 
\item {[Br3]} K. Brakke, Minimal surfaces, corners, and wires.
J. Geom. Anal. 2 (1992), no. 1, 11--36. 
\item {[Br4]} K. Brakke, Soap films and covering spaces.  
J. Geom. Anal.  5  (1995),  no. 4, 445--514.
\item {[DMS]} G. Dal Maso, J.-M. Morel, and S. Solimini, 
A variational method in image segmentation: Existence and
approximation results. 
Acta Math. 168 (1992), no. 1-2, 89--151.
\item {[Da1]} G. David, Limits of Almgren-quasiminimal sets. 
Proceedings of the conference on Harmonic Analysis, 
Mount Holyoke, A.M.S. Contemporary Mathematics series, Vol. 320 
(2003), 119--145.
\item {[Da2]} G. David, Singular sets of minimizers for 
the Mumford-Shah functional.
Progress in Mathematics 233 (581p.), Birkh\"auser 2005.
\item {[Da3]} G. David, Quasiminimal sets for Hausdorff measures. 
Recent developments in nonlinear partial differential equations, 81--99, 
Contemp. Math., 439, Amer. Math. Soc., Providence, RI, 2007.
\item {[Da4]} G. David, Low regularity for almost-minimal sets in $\R^3$. 
Annales de la Facult\'{e} des Sciences de Toulouse, 
Vol 18, 1 (2009), 65--246.
\item {[Da5]} G. David, $C^{1+\alpha}$-regularity
for two-dimensional almost-minimal sets in $\R^n$.
J. Geom. Anal. 20 (2010), no. 4, 837--954.
\item {[Da6]} G. David, 
Regularity of minimal and almost minimal sets and 
cones: J. Taylor's theorem for beginners.
Lecture notes for the S\'{e}minaires de Math\'{e}matiques 
Sup\'{e}rieures,
Espaces m\'{e}triques mesur\'{e}s : aspects g\'{e}om\'{e}triques et 
analytiques, held in Montreal in July  1011, to be published.
\item {[DDT]} G. David, T. De Pauw, and T. Toro,
A generalization of Reifenberg's Theorem in $\Bbb R^3$.
Geometric And Functional Analysis 18 (2008), 1168-1235.
\item {[DS1]} G. David and S. Semmes, \underbar{Analysis of and on 
uniformly rectifiable sets}. A.M.S. series of
Mathematical surveys and monographs, Volume 38, 1993. 
\item {[DS2]} G. David and S. Semmes, Uniform rectifiability and 
quasiminimizing sets of arbitrary codimension. 
Memoirs of the A.M.S. Number 687, volume 144,  2000.
\item {[Dp]} T. De Pauw, Size minimizing surfaces.
Ann. Sci. \'{E}c. Norm. Sup\'{e}r. (4) 42 (2009), no. 1, 37--101.
\item {[DpH]} T. De Pauw and R. Hardt, Size minimization and 
approximating problems.
Calc. Var. Partial Differential Equations 17 (2003), 405--442.
\item {[Do]} J. Douglas, Solutions of the Plateau problem. 
Trans. Amer. Math. Soc. 33 (1931), no. 1, 263--321. 
\item {[Fe1]} H. Federer, \underbar{Geometric measure theory}.
Grundlehren der Mathematishen Wissenschaften 
153, Springer Verlag 1969.
\item {[Fe2]} H. Federer,
The singular sets of area minimizing rectifiable currents 
with codimension one and of area minimizing flat chains 
modulo two with arbitrary codimension.
Bull. Amer. Math. Soc. 76 (1970) 767Ð771.
\item {[FF]} H. Federer and W.H. Fleming, Normal and integral currents.
Ann. of Math.(2) 72 (1960), 458--520.
\item {[Fl]} Wendell H. Fleming,
On the oriented Plateau problem.
Rend. Circ. Mat. Palermo (2) 11 (1962), 69Ð90.
\item {[Fv1]} V. Feuvrier, Un r\'{e}sultat d'existence pour les 
ensembles minimaux par optimisation sur des grilles poly\'{e}drales.
Th\`{e}se de l'universit\'{e} de Paris-Sud 11, Orsay, Septembre 2008.
\item {[Fv2]} V. Feuvrier, Remplissage de l'espace Euclidien par des 
complexes poly\'{e}driques d'orientation impos\'{e}e et de 
rotondit\'{e} uniforme. Preprint, 2008, arXiv:0812.4709.
\item {[Fv3]} V. Feuvrier, Condensation of polyhedric structures onto 
soap films. Preprint, 2009, arXiv:0906.3505.
\item {[Ga]} R. Garnier, Le probl\`{e}me de Plateau.
Annales Scientifiques de l'Ecole Normale Sup\'{e}rieure,
vol. 45 (1928), pp. 53--144.
\item {[HS]} R. Hardt and L. Simon, 
Boundary regularity and embedded solutions for the oriented Plateau problem. 
Ann. of Math. 110 (1979), 439--486.
\item {[Ha1]} J. Harrison,
Cartan's magic formula and soap film structures.
Journal of Geometric Analysis 14 (2004) no 1 47-61.
\item {[Ha2]} J. Harrison,
On Plateau's problem for soap films with a bound on energy.
Journal of Geometric Analysis 14 (2004) no 2, 319-329.
\item {[Ha3]} J. Harrison, Solution of Plateau's problem. Preprint.
\item {[He]} A. Heppes, Isogonal sph\"arischen netze. 
Ann. Univ. Sci. Budapest, E\"otv\"os Sect. Math. 7 (1964), 41--48.
\item {[Lam]} E. Lamarle, 
Sur la stabilit\'{e} des syst\`{e}mes liquides en lames minces. 
M\'{e}m. Acad. R. Belg. 35 (1864), 3--104.
\item {[Lawl]}  Gary Lawlor, 
Pairs of planes which are not size-minimizing. 
Indiana Univ. Math. J. 43 (1994), 651--661. 
\item {[LM1]}  Gary Lawlor and Frank Morgan, 
Paired calibrations applied to soap films, immiscible fluids, and 
surfaces or networks minimizing other norms.
Pacific J. Math.  166  (1994),  no. 1, 55--83.
\item {[LM2]} G. Lawlor and F. Morgan, 
Curvy slicing proves that triple junctions locally minimize area.
J. Diff. Geom. 44 (1996), 514--528. 
\item {[Laws]} H. Blaine Lawson Jr, 
Lectures on minimal submanifolds. Vol. I. Second edition. 
Mathematics Lecture Series, 9. 
Publish or Perish, Inc., Wilmington, Del., 1980. 
iv+178 pp. ISBN: 0-914098-18-7.
\item {[Li1]} X. Liang, 
Ensembles et c\^{o}nes minimaux de dimension $2$ 
dans les espaces euclidiens. 
Thesis, Universit\'{e} de Paris-Sud 11, Orsay, December 2010. 
\item {[Li2]} X. Liang, 
Almgren-minimality of unions of two almost orthogonal planes in $\R^4$.
Preprint, Universit\'{e} de Paris-Sud 11, Orsay, 2011 arXiv:1103.1468.
\item {[Li3]} X. Liang, 
Topological minimal sets and their applications.
Preprint, Universit\'{e} de Paris-Sud 11, 
Orsay, 2011, arXiv:1103.3871. 
\item {[Lu]} T. D. Luu, 
R\'{e}gularit\'{e} des c\^{o}nes et ensembles minimaux 
de dimension $3$ dans $\R^4$.
Thesis, Universit\'{e} de Paris 11, Orsay 2011. 
\item {[Ma]}  P. Mattila, \underbar{Geometry of sets and 
measures in Euclidean space}. Cambridge Studies in
Advanced Mathematics 44, Cambridge University Press l995.
\item {[Mo1]} F. Morgan, Size-minimizing rectifiable currents.
Invent. Math. 96 (1989), no. 2, 333--348.
\item {[Mo2]} F. Morgan, Minimal surfaces, crystals, shortest networks, and 
undergraduate research. Math. Intelligencer 14 (1992), no. 3, 37--44.
\item {[Mo3]} F. Morgan, 
Soap films and mathematics. 
Proceedings of Symposia in Pure Mathematics, 54, Part 1, (1993). 
\item {[Mo4]} F. Morgan, 
$(M,\varepsilon,\delta)$-minimal curve regularity. 
Proc. Amer. Math. Soc. 120  (1994),  no. 3, 677--686. 
\item {[Mo5]} F. Morgan, 
\underbar{Geometric measure theory. A beginner's guide}. 
4th edition, Academic Press, Inc., San Diego, CA, 2009. 264pp.
\item {[Ni1]} Johannes C. C. Nitsche,
Vorlesungen Ÿber MinimalflŠchen. 
Die Grundlehren der mathematischen Wissenschaften, Band 199. 
Springer-Verlag, Berlin-New York, 1975. xiii+775 pp. 
\item {[Ni2]} Johannes C. C. Nitsche,
Lectures on minimal surfaces.
Vol. 1. Introduction, fundamentals, geometry and 
basic boundary value problems, 
Translated from the German by Jerry M. Feinberg. 
With a German foreword. 
Cambridge University Press, Cambridge, 1989. 
xxvi+563 pp. ISBN: 0-521-24427-7.
\item {[Os]} R. Osserman, 
A survey of minimal surfaces, Second edition.
Dover Publications, Inc., New York, 1986. vi+207 pp. 
\item {[Pl]} J. Plateau, 
Statique exp\'{e}rimentale et th\'{e}orique des liquides 
soumis aux seules forces mol\'{e}culaires. Gauthier-Villars, Paris, 1873.
\item {[Ra1]} T. Rad\'{o}, The problem of least area and the problem 
of Plateau. Math. Z. 32 (1930), 763--796. 
\item {[Ra2]} T. Rad\'{o}, On the problem of Plateau. Ergebnisse der 
Mathematik und ihrer Grenzgebiete, Vol. 2, Springer, Berlin 1933; 
Reprinted Chelsea, New York 1951. 
\item {[Re1]} E. R. Reifenberg, 
Solution of the Plateau Problem for $m$-dimensional surfaces 
of varying topological type.
Acta Math. 104, 1960, 1--92.
\item {[Re2]} E. R. Reifenberg, 
An epiperimetric inequality related to the analyticity of minimal surfaces. 
Ann. of Math. (2) 80, 1964, 1--14.
\item {[Re3]} E. R. Reifenberg, On the analyticity of minimal surfaces.
Annals of Math. 80 (1964), 15--21.
\item {[Ta1]} J. Taylor, Regularity of the singular sets of 
two-dimensional area-minimizing flat chains modulo $3$ in $\R^{3}$.
Invent. Math. 22 (1973), 119--159.
\item {[Ta2]} J. Taylor, The structure of singularities in 
soap-bubble-like  and soap-film-like minimal surfaces.
Ann. of Math. (2) 103 (1976), no. 3, 489--539.
\item {[Wh]} B. White, 
Regularity of area-minimizing hypersurfaces at boundaries with 
multiplicity. In ``Seminar on Minimal Submanifolds", E. Bombieri Ed., 
pp. 293--301, Princeton University Press, NJ. 1983.

\ms\ms\noindent
Guy David,
\par\noindent
Univ Paris-Sud, Laboratoire de Math\'{e}matiques, UMR 8628,
Orsay, F-91405
\par\noindent
CNRS, Orsay, F-91405
\par\noindent
and Institut Universitaire de France.

\bye